\documentclass[a4paper,twoside,10pt]{amsart}
\usepackage{amsmath,amsfonts,amssymb,amsthm}
\newtheorem{theorem}{Theorem}[section]

\usepackage{mathrsfs}
\usepackage{color,graphicx}
\usepackage[a4paper]{geometry}
\numberwithin{equation}{section}

\author[G. Nemes]{Gerg\H{o} Nemes}
\address{Central European University, Department of Mathematics and its Applications, H-1051 Budapest, N\'ador utca 9, Hungary}
\email{nemesgery@gmail.com}

\keywords{asymptotic expansions, Anger--Weber function, error bounds, Stokes phenomenon, late coefficients.}
\subjclass[2010]{41A60, 30E15, 34M40}
\begin{document}

\title[Resurgence of the Anger--Weber function]{The resurgence properties of\\ the large order asymptotics of\\ the Anger--Weber function I}

\begin{abstract} The aim of this paper is to derive new representations for the Anger--Weber function, exploiting the reformulation of the method of steepest descents by C. J. Howls (Howls, Proc. R. Soc. Lond. A \textbf{439} (1992) 373--396). Using these representations, we obtain a number of properties of the large order asymptotic expansions of the Anger--Weber function, including explicit and realistic error bounds, asymptotics for the late coefficients, exponentially improved asymptotic expansions, and the smooth transition of the Stokes discontinuities.
\end{abstract}
\maketitle

\section{Introduction}\label{section1}

In this paper, we investigate the large $\nu$ asymptotics of the Anger--Weber function $\mathbf{A}_{ - \nu } \left( \nu x\right)$. The asymptotic expansion of this function has different forms according to whether $0<x<1$, $x=1$ or $x>1$ \cite[p. 298]{NIST} (see also Olver \cite[p. 352]{Olver}). We shall consider the latter two cases. A brief discussion about the expansion for $\mathbf{A}_{ - \nu } \left( \nu x\right)$ with $0<x<1$ is given in Section \ref{section6}. Meijer \cite{Meijer} gave error bounds for the asymptotic expansion of $\mathbf{A}_{ - \nu } \left( \nu x\right)$ when $x > 1$. Dingle \cite{Dingle} obtained exponentially improved versions of the asymptotic series and asymptotic approximations for their late terms. Nevertheless, the derivation of his results is based on
interpretive, rather than rigorous, methods. In an earlier paper \cite{Nemes}, we proved resurgence-type formulas for the Hankel function $H_{\nu}^{\left( 1 \right)} \left( \nu x\right)$ when $x=1$ and $x>1$, respectively. The main aim of this paper is to derive similar new representations for the Anger--Weber function $\mathbf{A}_{ - \nu } \left( \nu x\right)$. Our derivation is based on the reformulation of the method of steepest descents by Howls \cite{Howls}. Using these representations, we obtain a number of properties of the large order asymptotic expansions of the Anger--Weber function, including explicit and realistic error bounds, asymptotics for the late coefficients, exponentially improved asymptotic expansions, and the smooth transition of the Stokes discontinuities. Some of our error bounds coincide with the ones given by Meijer while some others are simpler. Our analysis also provides a rigorous treatment of Dingle's formal expansions.

Our first theorem describes the resurgence properties of the asymptotic expansion of $\mathbf{A}_{ - \nu } \left( \nu x\right)$ for $x > 1$. We employ the substitution $x = \sec \beta$ with an appropriate $0 < \beta <\frac{\pi}{2}$. The notations follow the ones given in \cite[p. 298]{NIST}. Throughout this paper, empty sums are taken to be zero.

\begin{theorem}\label{thm1} Let $0 < \beta <\frac{\pi}{2}$ be a fixed acute angle, and let $N$ be a non-negative integer. Then we have
\begin{equation}\label{eq18}
\mathbf{A}_{ - \nu } \left( {\nu \sec \beta } \right) =  - \frac{1}{\pi }\sum\limits_{n = 0}^{N - 1} {\frac{{\left( {2n} \right)!a_n \left( { - \sec \beta } \right)}}{{\nu ^{2n + 1} }}}  + R_N \left( {\nu ,\beta } \right)
\end{equation}
for $-\frac{\pi}{2} < \arg \nu < \frac{\pi}{2}$, with
\begin{equation}\label{eq9}
a_n \left( { - \sec \beta } \right) =  - \frac{1}{{\left( {2n} \right)!}}\left[ {\frac{{d^{2n} }}{{dt^{2n} }}\left( {\frac{t}{{\sec \beta \sinh t - t}}} \right)^{2n + 1} } \right]_{t = 0} = \frac{{\left( { - 1} \right)^{n + 1} }}{{\left( {2n} \right)!}}\int_0^{ + \infty } {t^{2n} iH_{it}^{\left( 1 \right)} \left( {it\sec \beta } \right)dt} 
\end{equation}
and
\begin{equation}\label{eq8}
R_N \left( {\nu ,\beta } \right) = \frac{{\left( { - 1} \right)^N }}{{\pi \nu ^{2N + 1} }}\int_0^{ + \infty } {\frac{{t^{2N} }}{{1 + \left( {t/\nu } \right)^2 }}iH_{it}^{\left( 1 \right)} \left( {it\sec \beta } \right)dt} .
\end{equation}
\end{theorem}

It was shown in \cite{Nemes}, that for any fixed $0 < \beta <\frac{\pi}{2}$ and non-negative integer $M$, the Hankel function $H_\nu ^{\left( 1 \right)} \left( {\nu \sec \beta } \right)$ has the representation
\begin{equation}\label{eq40}
H_\nu ^{\left( 1 \right)} \left( {\nu \sec \beta } \right) = \frac{{e^{i\nu \left( {\tan \beta  - \beta } \right) - \frac{\pi}{4}i} }}{{\left( {\frac{1}{2}\nu \pi \tan \beta } \right)^{\frac{1}{2}} }}\left( {\sum\limits_{m = 0}^{M - 1} {\left( { - 1} \right)^m \frac{{U_m \left( {i\cot \beta } \right)}}{{\nu ^m }}}  + R_M^{\left(H\right)} \left( {\nu ,\beta } \right)} \right),
\end{equation}
for $-\frac{\pi}{2} < \arg \nu < \frac{3\pi}{2}$ with
\begin{gather}\label{eq41}
\begin{split}
U_m \left( {i\cot \beta } \right) & = \left( { - 1} \right)^m \frac{{\left( {i\cot \beta } \right)^m }}{{2^m m!}}\left[ {\frac{{d^{2m} }}{{dt^{2m} }}\left( {\frac{1}{2}\frac{{t^2 }}{{i\cot \beta \left( {t - \sinh t} \right) + \cosh t - 1}}} \right)^{m + \frac{1}{2}} } \right]_{t = 0}\\
& = \frac{i^m}{2\left( {2\pi \cot \beta } \right)^{\frac{1}{2}}}\int_0^{ + \infty } {t^{m - \frac{1}{2}} e^{ -t \left( {\tan \beta  - \beta } \right)} \left( {1 + e^{ - 2\pi t} } \right)i H_{it}^{\left( 1 \right)} \left( {it\sec \beta } \right)dt} .
\end{split}
\end{gather}
The remainder term $R_M^{\left(H\right)} \left( {\nu ,\beta } \right)$ can be expressed as
\begin{equation}\label{eq42}
R_M^{\left(H\right)} \left( {\nu ,\beta } \right) =  \frac{1}{{2\left( {2\pi \cot \beta } \right)^{\frac{1}{2}} \left(i\nu\right)^M }}\int_0^{ + \infty } {\frac{{t^{M - \frac{1}{2}} e^{ - t\left( {\tan \beta  - \beta } \right)} }}{{1 + it/\nu }}\left( {1 + e^{ - 2\pi t} } \right)i H_{it}^{\left( 1 \right)} \left( {it\sec \beta } \right)dt} .
\end{equation}
This representation of the Hankel function will play an important role in later sections of this paper.

If $\mathbf{J}_{\nu}\left(z\right)$ denotes the Anger function, then $\mathbf{J}_{-\nu}\left(z\right) = \mathbf{J}_{\nu}\left(-z\right)$ and $\sin \left( {\pi \nu } \right)\mathbf{A}_{ - \nu } \left( {\nu \sec \beta } \right) = J_{ - \nu } \left( {\nu \sec \beta } \right) - \mathbf{J}_{ - \nu } \left( {\nu \sec \beta } \right)$ (see \cite[p. 296]{NIST}). From these and the continuation formulas for the Bessel and Hankel functions (see \cite[p. 222 and p. 226]{NIST}), we find
\begin{align*}
\sin \left( {\pi \nu } \right)\mathbf{A}_{ - \nu } \left( {\nu e^{2\pi im} \sec \beta } \right) = &\; J_{ - \nu } \left( {\nu e^{2\pi im} \sec \beta } \right) - \mathbf{J}_{ - \nu } \left( {\nu e^{2\pi im} \sec \beta } \right)\\
= &\; \sin \left( {\pi \nu } \right)\mathbf{A}_{ - \nu } \left( {\nu \sec \beta } \right) + \left( {e^{ - 2\pi im\nu }  - 1} \right)J_{ - \nu } \left( {\nu \sec \beta } \right)\\
= &\; \sin \left( {\pi \nu } \right)\mathbf{A}_{ - \nu } \left( {\nu \sec \beta } \right) - ie^{ - \pi i\left( {m - 1} \right)\nu } \sin \left( {\pi m\nu } \right)H_\nu ^{\left( 1 \right)} \left( {\nu \sec \beta } \right)\\ & - ie^{ - \pi i\left( {m + 1} \right)\nu } \sin \left( {\pi m\nu } \right)H_\nu ^{\left( 2 \right)} \left( {\nu \sec \beta } \right)
\end{align*}
for every integer $m$. From this expression and the resurgence formulas \eqref{eq18}, \eqref{eq40}, we can derive analogous representations in sectors of the form
\[
\left( {2m - \frac{1}{2}} \right)\pi  < \arg \nu  < \left( {2m + \frac{1}{2}} \right)\pi ,\; m \in \mathbb{Z}.
\]

Similarly, applying the continuation formulas
\begin{align*}
& - \sin \left( {\pi \nu } \right)\mathbf{A}_{\nu } \left( {\nu e^{\left( {2m + 1} \right)\pi i} \sec \beta } \right) = J_{ \nu } \left( {\nu e^{\left( {2m + 1} \right)\pi i} \sec \beta } \right) - \mathbf{J}_{ \nu } \left( {\nu e^{\left( {2m + 1} \right)\pi i} \sec \beta } \right)\\ & = \sin \left( {\pi \nu } \right)\mathbf{A}_{ - \nu } \left( {\nu \sec \beta } \right) + e^{\left( {2m + 1} \right)\pi i\nu } J_\nu  \left( {\nu \sec \beta } \right) - J_{ - \nu } \left( {\nu \sec \beta } \right)\\
& = \sin \left( {\pi \nu } \right)\mathbf{A}_{ - \nu } \left( {\nu \sec \beta } \right) + ie^{\pi i\left( {m + 1} \right)\nu } \sin \left( {\pi m\nu } \right)H_\nu ^{\left( 1 \right)} \left( {\nu \sec \beta } \right) + ie^{\pi im\nu } \sin \left( {\pi \left( {m + 1} \right)\nu } \right)H_\nu ^{\left( 2 \right)} \left( {\nu \sec \beta } \right)
\end{align*}
and the representations \eqref{eq18}, \eqref{eq40}, we can obtain resurgence formulas in any sector of the form
\[
\left( {2m + \frac{1}{2}} \right)\pi  < \arg \nu  < \left( {2m + \frac{3}{2}} \right)\pi ,\; m \in \mathbb{Z}.
\]
The lines $\arg \nu  = \left( {2m \pm \frac{1}{2}} \right)\pi$ are the Stokes lines for the function $\mathbf{A}_{ - \nu } \left( {\nu \sec \beta } \right)$.

When $\nu$ is an integer, the limiting values have to be taken in these continuation formulas.

The second theorem provides a resurgence formula for $\mathbf{A}_{ - \nu } \left( \nu \right)$.

\begin{theorem}\label{thm2} For any non-negative integer $N$, we have
\begin{equation}\label{eq19}
\mathbf{A}_{ - \nu } \left( \nu  \right) = \frac{1}{{3\pi }}\sum\limits_{n = 0}^{N - 1} {d_{2n} \frac{{\Gamma \left( {\frac{{2n + 1}}{3}} \right)}}{{\nu ^{\frac{{2n + 1}}{3}} }}}  + R_N \left( \nu  \right)
\end{equation}
for $-\frac{3\pi}{2} < \arg \nu < \frac{3\pi}{2}$, with
\begin{equation}\label{eq17}
d_{2n} = \frac{1}{{\left( {2n} \right)!}}\left[ {\frac{{d^{2n} }}{{dt^{2n} }}\left( {\frac{{t^3 }}{{\sinh t - t}}} \right)^{\frac{{2n + 1}}{3}} } \right]_{t = 0}  = \frac{{\left( { - 1} \right)^n }}{{\Gamma \left( {\frac{{2n + 1}}{3}} \right)}}\int_0^{ + \infty } {t^{\frac{{2n - 2}}{3}} e^{ - 2\pi t} i H_{it}^{\left( 1 \right)} \left( {it} \right)dt}
\end{equation}
and
\begin{equation}\label{eq16}
R_N \left( \nu  \right) = \frac{{\left( { - 1} \right)^N }}{{3\pi \nu ^{\frac{{2N + 1}}{3}} }}\int_0^{ + \infty } {\frac{{t^{\frac{{2N - 2}}{3}} e^{ - 2\pi t} }}{{1 + \left( {t/\nu } \right)^{\frac{2}{3}} }}iH_{it}^{\left( 1 \right)} \left( {it} \right)dt} .
\end{equation}
The cube roots are defined to be positive on the positive real line and are defined by analytic continuation elsewhere.
\end{theorem}

In the previous paper \cite{Nemes}, we proved a similar representation for the Hankel function $H_\nu ^{\left( 1 \right)} \left( \nu  \right)$, in particular for any non-negative integer $N$, we have
\begin{equation}\label{eq53}
H_\nu ^{\left( 1 \right)} \left( \nu  \right) =  - \frac{2}{{3\pi }}\sum\limits_{n = 0}^{N - 1} {d_{2n} e^{\frac{{2\left( {2n + 1} \right)\pi i}}{3}} \sin \left( {\frac{{\left( {2n + 1} \right)\pi }}{3}} \right)\frac{{\Gamma \left( {\frac{{2n + 1}}{3}} \right)}}{{\nu ^{\frac{{2n + 1}}{3}} }}}  + R_N^{\left( H \right)} \left( \nu  \right)
\end{equation}
with $-\frac{\pi}{2} < \arg \nu < \frac{3\pi}{2}$. The remainder term $R_N^{\left( H \right)} \left( \nu  \right)$ has the integral representation
\begin{equation}\label{eq54}
R_N^{\left( H \right)} \left( \nu  \right) = \frac{{\left( { - 1} \right)^N }}{{3\pi \nu ^{\frac{{2N + 1}}{3}} }}\int_0^{ + \infty } {t^{\frac{{2N - 2}}{3}} e^{ - 2\pi t} \left( {\frac{{e^{\frac{{\left( {2N + 1} \right)\pi i}}{3}} }}{{1 + \left( {t/\nu } \right)^{\frac{2}{3}} e^{\frac{{2\pi i}}{3}} }} + \frac{1}{{1 + \left( {t/\nu } \right)^{\frac{2}{3}} }}} \right)H_{it}^{\left( 1 \right)} \left( {it} \right)dt} .
\end{equation}
The cube roots are defined to be positive on the positive real line and are defined by analytic continuation elsewhere. This result will be important for us in later sections of the paper.

Again, the formula \eqref{eq19} can be extended to other sectors of the complex plane. (One has to replace the factor $\sec \beta$ by 1 in the continuation formulas given earlier.)

If we neglect the remainder terms and extend the sums to $N = \infty$ in Theorems \ref{thm1} and \ref{thm2}, we recover the known asymptotic series of the Anger--Weber function. Some other formulas for the coefficients $a_n \left( { - \sec \beta } \right)$ can be found in Appendix \ref{appendixa}. For the computation of the $d_{2n}$, see \cite[Appendix A]{Nemes}.

In the following two theorems, we give exponentially improved asymptotic expansions for the function $\mathbf{A}_{ - \nu } \left(\nu x\right)$ when $x>1$ and $x=1$, respectively. These new expansions can be viewed as the mathematically rigorous forms of the terminated series of Dingle \cite[pp. 485]{Dingle}. We express these expansions in terms of the Terminant function $\widehat T_p\left(w\right)$ whose definition and basic properties are given in Section \ref{section5}. In Theorem \ref{thm3}, $R_N \left( {\nu ,\beta } \right)$ is defined by \eqref{eq18} and it is extended to the sector $\left|\arg \nu\right| \leq \frac{3\pi}{2}$ via analytic continuation. Throughout this paper, we use subscripts in the $\mathcal{O}$ notations to indicate the dependence of the implied constant on certain parameters.

\begin{theorem}\label{thm3} Suppose that $\left|\arg \nu\right| \leq \frac{3\pi}{2}$, $\left|\nu\right|$ is large and $N = \frac{1}{2}\left| \nu  \right|\left( {\tan \beta  - \beta } \right) + \rho$ is a positive integer with $\rho$ being bounded. Then
\begin{align*}
R_N \left( {\nu ,\beta } \right) = \; & i\frac{{e^{i\nu \left( {\tan \beta  - \beta } \right) - \frac{\pi }{4}i} }}{{\left( {\frac{1}{2}\nu \pi \tan \beta } \right)^{\frac{1}{2}} }}\sum\limits_{m = 0}^{M - 1} {\left( { - 1} \right)^m \frac{{U_m \left( {i\cot \beta } \right)}}{{\nu ^m }}\widehat T_{2N - m + \frac{1}{2}} \left( {i\nu \left( {\tan \beta  - \beta } \right)} \right)}\\ & - i\frac{{e^{ - i\nu \left( {\tan \beta  - \beta } \right) + \frac{\pi }{4}i} }}{{\left( {\frac{1}{2}\nu \pi \tan \beta } \right)^{\frac{1}{2}} }}\sum\limits_{m = 0}^{M - 1} {\frac{{U_m \left( {i\cot \beta } \right)}}{{\nu ^m }}\widehat T_{2N - m + \frac{1}{2}} \left( { - i\nu \left( {\tan \beta  - \beta } \right)} \right)}  + R_{N,M} \left( {\nu ,\beta } \right)
\end{align*}
with $M$ being an arbitrary fixed non-negative integer, and
\[
R_{N,M} \left( {\nu ,\beta } \right) = \mathcal{O}_{M,\rho } \left( {\frac{{e^{ - \left| \nu  \right|\left( {\tan \beta  - \beta } \right)} }}{{\left( {\frac{1}{2}\left| \nu  \right|\pi \tan \beta } \right)^{\frac{1}{2}} }}\frac{{\left| {U_M \left( {i\cot \beta } \right)} \right|}}{{\left| \nu  \right|^M }}} \right)
\]
for $\left|\arg \nu\right| \leq \frac{\pi}{2}$;
\[
R_{N,M} \left( {\nu ,\beta } \right) = \mathcal{O}_{M,\rho } \left( {\frac{{e^{ \mp \Im \left( \nu  \right)\left( {\tan \beta  - \beta } \right)} }}{{\left( {\frac{1}{2}\left| \nu  \right|\pi \tan \beta } \right)^{\frac{1}{2}} }}\frac{{\left| {U_M \left( {i\cot \beta } \right)} \right|}}{{\left| \nu  \right|^M }}} \right)
\]
for $\frac{\pi}{2} \leq \pm \arg \nu \leq \frac{3\pi}{2}$.
\end{theorem}

\begin{theorem}\label{thm4} Define $R_{N,M,K}\left(\nu\right)$ by
\begin{align*}
\mathbf{A}_{ - \nu } \left( \nu  \right) = \frac{1}{{3\pi \nu ^{\frac{1}{3}} }}\sum\limits_{n = 0}^{N - 1} {d_{6n} \frac{{\Gamma \left( {2n + \frac{1}{3}} \right)}}{{\nu ^{2n} }}} & + \frac{1}{{3\pi \nu }}\sum\limits_{m = 0}^{M - 1} {d_{6m + 2} \frac{{\Gamma \left( {2m + 1} \right)}}{{\nu ^{2m} }}} \\ & + \frac{1}{{3\pi \nu ^{\frac{5}{3}} }}\sum\limits_{k = 0}^{K - 1} {d_{6k + 4} \frac{{\Gamma \left( {2k + \frac{5}{3}} \right)}}{{\nu ^{2k} }}}  + R_{N,M,K} \left( \nu  \right),
\end{align*}
where
\[
N= \pi\left| \nu  \right| + \rho, \; M= \pi\left| \nu  \right| + \sigma  \; \text{ and } \; K= \pi\left| \nu  \right| + \eta,
\]
$\left|\nu\right|$ being large, $\rho$, $\sigma$ and $\eta$ being bounded quantities such that $N,M,K \geq 1$. Then
\begin{gather}\label{eq68}
\begin{split}
R_{N,M,K} \left( \nu  \right) =\; & i\frac{{e^{ - 2\pi i\nu } }}{3}\frac{2}{{3\pi }}\sum\limits_{j = 0}^{J - 1} {d_{2j} \sin \left( {\frac{{\left( {2j + 1} \right)\pi }}{3}} \right)\frac{{\Gamma \left( {\frac{{2j + 1}}{3}} \right)}}{{\nu ^{\frac{{2j + 1}}{3}} }}\widehat T_{2N - \frac{{2j}}{3}} \left( { - 2\pi i\nu } \right)} \\
& - ie^{\frac{\pi }{3}i} \frac{{e^{2\pi i\nu } }}{3}\frac{2}{{3\pi }}\sum\limits_{j = 0}^{J - 1} {d_{2j} e^{\frac{{2\left( {2j + 1} \right)\pi i}}{3}} \sin \left( {\frac{{\left( {2j + 1} \right)\pi }}{3}} \right)\frac{{\Gamma \left( {\frac{{2j + 1}}{3}} \right)}}{{\nu ^{\frac{{2j + 1}}{3}} }}\widehat T_{2N - \frac{{2j}}{3}} \left( {2\pi i\nu } \right)} \\
& + i\frac{{e^{ - 2\pi i\nu } }}{3}\frac{2}{{3\pi }}\sum\limits_{\ell  = 0}^{L - 1} {d_{2\ell } \sin \left( {\frac{{\left( {2\ell  + 1} \right)\pi }}{3}} \right)\frac{{\Gamma \left( {\frac{{2\ell  + 1}}{3}} \right)}}{{\nu ^{\frac{{2\ell  + 1}}{3}} }}\widehat T_{2M - \frac{{2\ell  - 2}}{3}} \left( { - 2\pi i\nu } \right)} \\
& + i\frac{{e^{2\pi i\nu } }}{3}\frac{2}{{3\pi }}\sum\limits_{\ell  = 0}^{L - 1} {d_{2\ell } e^{\frac{{2\left( {2\ell  + 1} \right)\pi i}}{3}} \sin \left( {\frac{{\left( {2\ell  + 1} \right)\pi }}{3}} \right)\frac{{\Gamma \left( {\frac{{2\ell  + 1}}{3}} \right)}}{{\nu ^{\frac{{2\ell  + 1}}{3}} }}\widehat T_{2M - \frac{{2\ell  - 2}}{3}} \left( {2\pi i\nu } \right)} \\
& + i\frac{{e^{ - 2\pi i\nu } }}{3}\frac{2}{{3\pi }}\sum\limits_{q = 0}^{Q - 1} {d_{2q} \sin \left( {\frac{{\left( {2q + 1} \right)\pi }}{3}} \right)\frac{{\Gamma \left( {\frac{{2q + 1}}{3}} \right)}}{{\nu ^{\frac{{2q + 1}}{3}} }}\widehat T_{2K - \frac{{2q - 4}}{3}} \left( { - 2\pi i\nu } \right)} \\
& - ie^{ - \frac{\pi }{3}i} \frac{{e^{2\pi i\nu } }}{3}\frac{2}{{3\pi }}\sum\limits_{q = 0}^{Q - 1} {d_{2q} e^{\frac{{2\left( {2q + 1} \right)\pi i}}{3}} \sin \left( {\frac{{\left( {2q + 1} \right)\pi }}{3}} \right)\frac{{\Gamma \left( {\frac{{2q + 1}}{3}} \right)}}{{\nu ^{\frac{{2q + 1}}{3}} }}\widehat T_{2K - \frac{{2q - 4}}{3}} \left( {2\pi i\nu } \right)} \\
& + R_{N,M,K}^{J,L,Q} \left( \nu  \right),
\end{split}
\end{gather}
where $J$, $L$ and $Q$ are arbitrary fixed non-negative integers satisfying $J,L,Q \equiv 0 \mod 3$, and
\begin{gather}\label{eq62}
\begin{split}
R_{N,M,K}^{J,L,Q} \left( \nu  \right) = \mathcal{O}_{J,\rho } \left( {e^{ - 2\pi \left| \nu  \right|} \left| {d_{2J} } \right|\frac{{\Gamma \left( {\frac{{2J + 1}}{3}} \right)}}{{\left| \nu  \right|^{\frac{{2J + 1}}{3}} }}} \right) & + \mathcal{O}_{L,\sigma } \left( {e^{ - 2\pi \left| \nu  \right|} \left| {d_{2L} } \right|\frac{{\Gamma \left( {\frac{{2L + 1}}{3}} \right)}}{{\left| \nu  \right|^{\frac{{2L + 1}}{3}} }}} \right)\\
& + \mathcal{O}_{Q,\eta } \left( {e^{ - 2\pi \left| \nu  \right|} \left| {d_{2Q} } \right|\frac{{\Gamma \left( {\frac{{2Q + 1}}{3}} \right)}}{{\left| \nu  \right|^{\frac{{2Q + 1}}{3}} }}} \right)
\end{split}
\end{gather}
for $-\frac{\pi}{2} \leq \arg\nu \leq \frac{\pi}{2}$;
\begin{gather}\label{eq63}
\begin{split}
R_{N,M,K}^{J,L,Q} \left( \nu  \right) = \mathcal{O}_{J,\rho } \left( {e^{ \mp 2\pi \Im \left( \nu  \right)} \left| {d_{2J} } \right|\frac{{\Gamma \left( {\frac{{2J + 1}}{3}} \right)}}{{\left| \nu  \right|^{\frac{{2J + 1}}{3}} }}} \right) & + \mathcal{O}_{L,\sigma } \left( {e^{ \mp 2\pi \Im \left( \nu  \right)} \left| {d_{2L} } \right|\frac{{\Gamma \left( {\frac{{2L + 1}}{3}} \right)}}{{\left| \nu  \right|^{\frac{{2L + 1}}{3}} }}} \right)\\
& + \mathcal{O}_{Q,\eta } \left( {e^{ \mp 2\pi \Im \left( \nu  \right)} \left| {d_{2Q} } \right|\frac{{\Gamma \left( {\frac{{2Q + 1}}{3}} \right)}}{{\left| \nu  \right|^{\frac{{2Q + 1}}{3}} }}} \right)
\end{split}
\end{gather}
for $\frac{\pi}{2} \leq \pm \arg \nu \leq \frac{3\pi}{2}$;
\begin{align*}
R_{N,M,K}^{J,L,Q} \left( \nu  \right) = \; & \mathcal{O}_{J,\rho } \left( {\cosh \left( {2\pi \Im \left( \nu  \right)} \right)\left| {d_{2J} } \right|\frac{{\Gamma \left( {\frac{{2J + 1}}{3}} \right)}}{{\left| \nu  \right|^{\frac{{2J + 1}}{3}} }}} \right) + \mathcal{O}_{L,\sigma } \left( {\cosh \left( {2\pi \Im \left( \nu  \right)} \right)\left| {d_{2L} } \right|\frac{{\Gamma \left( {\frac{{2L + 1}}{3}} \right)}}{{\left| \nu  \right|^{\frac{{2L + 1}}{3}} }}} \right)\\
& +\mathcal{O}_{Q,\eta } \left( {\cosh \left( {2\pi \Im \left( \nu  \right)} \right)\left| {d_{2Q} } \right|\frac{{\Gamma \left( {\frac{{2Q + 1}}{3}} \right)}}{{\left| \nu  \right|^{\frac{{2Q + 1}}{3}} }}} \right) + \mathcal{O}_{J,L,Q} \left( {\left| \nu  \right|^{ - \frac{1}{3}} } \right)
\end{align*}
for $\frac{3\pi}{2} \leq \pm \arg \nu \leq \frac{5\pi}{2}$. Moreover, if $J=L=Q$, then the bound \eqref{eq62} remains valid in the larger sector $-\frac{3\pi}{2} \leq  \arg \nu \leq \frac{3\pi}{2}$, and the estimate \eqref{eq63} holds in the sectors $\frac{3\pi}{2} \leq \mp \arg \nu \leq \frac{5\pi}{2}$.
\end{theorem}

The assumption that $J,L,Q \equiv 0 \mod 3$ is only for simplicity. Estimations for $R_{N,M,K}^{J,L,Q} \left( \nu  \right)$ when $J$, $L$ or $Q$ may not be divisible by $3$ can be obtained similarly.

We remark that Dingle writes $A_\nu\left(z\right)$ in place of $\mathbf{A}_{-\nu}\left(z\right)$; and Olver's definition for $\mathbf{A}_{-\nu}\left(z\right)$ omits the factor $\frac{1}{\pi}$ in \eqref{eq69} below.

The rest of the paper is organized as follows. In Section \ref{section2}, we prove the resurgence formulas stated in Theorems \ref{thm1} and \ref{thm2}. In Section \ref{section3}, we give explicit and realistic error bounds for the asymptotic expansions of $\mathbf{A}_{-\nu}\left(\nu x\right)$ when $x\geq 1$ using the results of Section \ref{section2}. In Section \ref{section4}, asymptotic approximations for $a_n \left( { - \sec \beta } \right)$ as $n \to +\infty$ are given. In Section \ref{section5}, we prove the exponentially improved expansions presented in Theorems \ref{thm3} and \ref{thm4}, and provide a detailed discussion of the Stokes phenomenon related to the expansions of $\mathbf{A}_{-\nu}\left( \nu x \right)$. The paper concludes with a discussion in Section \ref{section6}.

\section{Proofs of the resurgence formulas}\label{section2}

Our analysis is based on the integral definition of the Anger--Weber function
\begin{equation}\label{eq69}
\mathbf{A}_{ - \nu } \left( z \right) = \frac{1}{\pi }\int_0^{ + \infty } {e^{\nu t - z\sinh t} dt} \quad \left| {\arg z} \right| < \frac{\pi }{2}.
\end{equation}
If $z = \nu x$, where $x$ is a positive constant, then
\begin{equation}\label{eq10}
\mathbf{A}_{ - \nu } \left( {\nu x} \right) = \frac{1}{\pi }\int_0^{ + \infty } {e^{ - \nu \left( {x\sinh t - t} \right)} dt} \quad \left| {\arg \nu } \right| < \frac{\pi }{2}.
\end{equation}
The analysis is significantly different according to whether $x > 1$ or $x = 1$. The saddle points of the integrand are the roots of the equation $x \cosh t = 1$. Hence, the saddle points are given by $t_{\pm}^{\left( k \right)} = \pm \mathrm{sech}^{ - 1} x + 2\pi ik$ where $k$ is an arbitrary integer. When $x=1$, we shall use the simpler notation $t^{\left( k \right)} = 2\pi ik$. We denote by $\mathscr{C}_{\pm}^{\left( k \right)}\left(\theta\right)$ the portion of the steepest paths that pass through the saddle point $t_{\pm}^{\left( k \right)}$. Here, and subsequently, we write $\theta = \arg \nu$. Similarly, $\mathscr{C}^{\left( k \right)}\left(\theta\right)$ denotes the steepest paths through the saddle point $t^{\left( k \right)}$. As for the path of integration in \eqref{eq10}, we take
\begin{equation}\label{eq71}
\mathscr{P}\left(\theta \right) = \left\{ {t \in \mathbb{C}:\arg \left[ {e^{i\theta } \left( {x\sinh t - t} \right)} \right] = 0,\, \Re \left( t \right) > 0,\, \left| {\Im \left( t \right)} \right| < \frac{\pi }{2}} \right\} .
\end{equation}
We remark that $\mathscr{P}\left(0\right)$ is the positive real axis. If $x=1$, the path $\mathscr{P}\left(\theta\right)$ is part of the contour $\mathscr{C}^{\left( 0 \right)}\left(\theta\right)$.

\subsection{Case (i): $x>1$} 
\begin{figure}[!t]
\def\svgwidth{0.6\columnwidth}
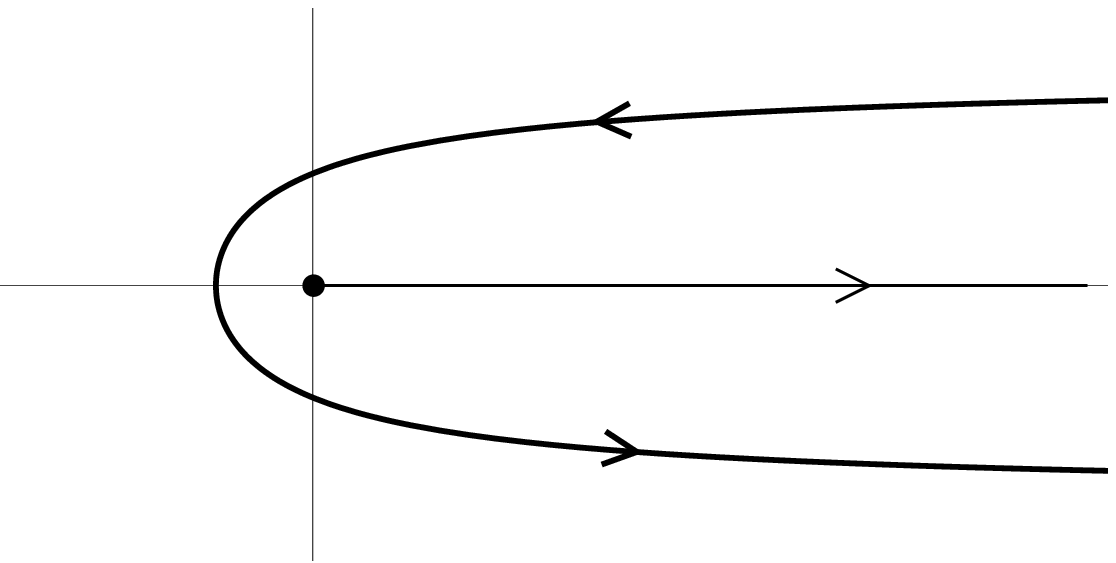
\caption{The contour $\Gamma$ encircling the path $\mathscr{P}\left(0\right)$.}
\label{fig1}
\end{figure}

Let $0<\beta<\frac{\pi}{2}$ be defined by $\sec \beta = x$. For simplicity, we assume that $\theta = 0$. In due course, we shall appeal to an analytic continuation argument to extend our results to complex $\nu$. Let $f\left( {t,\beta } \right) = \sec \beta \sinh t - t$. If
\begin{equation}\label{eq1}
\tau  = f\left( {t,\beta } \right),
\end{equation}
then $\tau$ is real on the curve $\mathscr{P}\left(0\right)$, and, as $t$ travels along this curve from $0$ to $+\infty $, $\tau$ increases from $0$ to $+\infty$. Therefore, corresponding to each positive value of $\tau$, there is a value of $t$, say $t\left(\tau\right)$, satisfying \eqref{eq1} with $t\left(\tau\right)>0$. In terms of $\tau$, we have
\[
\mathbf{A}_{ - \nu } \left( {\nu \sec \beta } \right) = \frac{1}{\pi }\int_0^{ + \infty } {e^{ - \nu \tau } \frac{{dt\left( \tau  \right)}}{{d\tau }}d\tau }  = \frac{1}{\pi }\int_0^{ + \infty } {e^{ - \nu \tau } \frac{1}{{\sec \beta \cosh t\left( \tau  \right) - 1}}d\tau } .
\]
Following Howls, we express the function involving $t\left(\tau\right)$ as a contour integral using the residue theorem, to find
\[
\mathbf{A}_{ - \nu } \left( {\nu \sec \beta } \right) = \frac{1}{\pi }\int_0^{ + \infty } {e^{ - \nu \tau } \frac{1}{{2\pi i}}\oint_\Gamma  {\frac{{f^{ - 1} \left( {u,\beta } \right)}}{{1 - \tau ^2 f^{ - 2} \left( {u,\beta } \right)}}du} d\tau }
\]
where the contour $\Gamma$ encircles the path $\mathscr{P}\left(0\right)$ in the positive direction and does not enclose any of the saddle points $t_ \pm ^{\left( k \right)}$ (see Figure \ref{fig1}). Now, we employ the well-known expression for non-negative integer $N$
\begin{equation}\label{eq2}
\frac{1}{1 - z} = \sum\limits_{n = 0}^{N-1} {z^n}  + \frac{z^N}{1 - z},\; z \neq 1,
\end{equation}
to expand the function under the contour integral in powers of $\tau ^2 f^{ - 2} \left( {u,\beta } \right)$. The result is
\[
\mathbf{A}_{ - \nu } \left( {\nu \sec \beta } \right) =  - \frac{1}{\pi}\sum\limits_{n = 0}^{N - 1} {\int_0^{ + \infty } {\tau^{2n} e^{ - \nu \tau }\frac{{ - 1}}{{2\pi i}}\oint_\Gamma  {\frac{{du}}{{f^{2n + 1} \left( {u,\beta } \right)}}} d\tau } }  + R_N \left( {\nu ,\beta } \right),
\]
where
\begin{equation}\label{eq3}
R_N \left( {\nu ,\beta } \right) = \frac{1}{\pi}\int_0^{ + \infty } {\tau^{2N} e^{ - \nu \tau } \frac{1}{2\pi i}\oint_\Gamma  {\frac{{f^{ - 2N - 1} \left( {u,\beta } \right)}}{{1 - \tau ^2 f^{ - 2} \left( {u,\beta } \right)}}du} d\tau } .
\end{equation}
The path $\Gamma$ in the sum can be shrunk into a small circle around $0$, and we arrive at
\begin{equation}\label{eq4}
\mathbf{A}_{ - \nu } \left( {\nu \sec \beta } \right) =  - \frac{1}{\pi }\sum\limits_{n = 0}^{N - 1} {\frac{{\left( {2n} \right)!a_n \left( { - \sec \beta } \right)}}{{\nu ^{2n + 1} }}}  + R_N \left( {\nu ,\beta } \right),
\end{equation}
where
\[
a_n \left( { - \sec \beta } \right) = \frac{{ - 1}}{{2\pi i}}\oint_{\left( {0^ +  } \right)} {\frac{{du}}{{f^{2n + 1} \left( {u,\beta } \right)}}}  =  - \frac{1}{{\left( {2n} \right)!}}\left[ {\frac{{d^{2n} }}{{dt^{2n} }}\left( {\frac{t}{{\sec \beta \sinh t - t}}} \right)^{2n + 1} } \right]_{t = 0} .
\]

\begin{figure}[t]
\def\svgwidth{0.7\columnwidth}
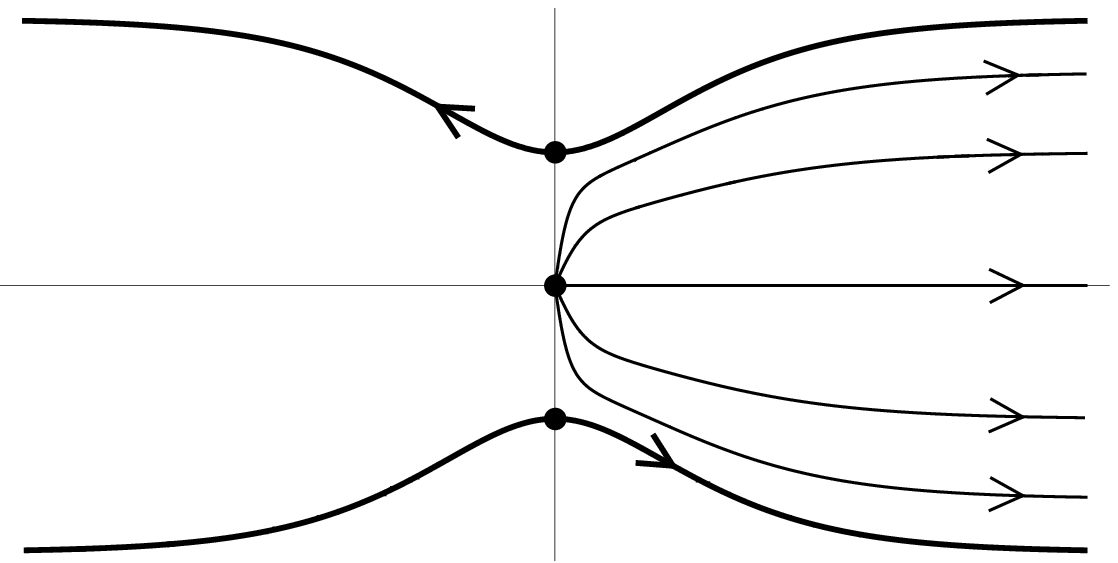
\caption{The path $\mathscr{P}\left(\theta\right)$ emanating from the origin when (i) $\theta=0$, (ii) $\theta=-\frac{\pi}{4}$, (iii) $\theta=-\frac{2\pi}{5}$, (iv) $\theta=\frac{\pi}{4}$, (v) $\theta=\frac{2\pi}{5}$. The paths $\mathscr{C}_+ ^{\left( 0 \right)} \left( { \frac{\pi }{2}} \right)$ and $\mathscr{C}_ - ^{\left( 0 \right)} \left( { - \frac{\pi }{2}} \right)$ are the adjacent contours to $0$. The domain $\Delta$ comprises all points between these two paths in the right-half plane.}
\label{fig2}
\end{figure}

Performing the change of variable $\nu \tau = s$ in \eqref{eq3} yields
\begin{equation}\label{eq5}
R_N \left( {\nu ,\beta } \right) = \frac{1}{{\pi \nu ^{2N + 1} }}\int_0^{ + \infty } {s^{2N} e^{ - s} \frac{1}{{2\pi i}}\oint_\Gamma  {\frac{{f^{ - 2N - 1} \left( {u,\beta } \right)}}{{1 - \left( {s/\nu } \right)^2 f^{ - 2} \left( {u,\beta } \right)}}du} ds} .
\end{equation}
This representation of $R_N \left( {\nu ,\beta } \right)$ and the formula \eqref{eq4} can be continued analytically if we choose $\Gamma = \Gamma\left(\theta\right)$ to be an infinite contour that surrounds the path $\mathscr{P}\left(\theta\right)$ in the anti-clockwise direction and that does not encircle any of the saddle points $t_ \pm ^{\left( k \right)}$. This continuation argument works until the path $\mathscr{P}\left(\theta\right)$ runs into a saddle point. In the terminology of Howls, such saddle points are called adjacent to the endpoint $0$. As
\[
\left|\arg \left( {f\left(0 ,\beta\right) - f\left( {t_ \pm ^{\left( k \right)} ,\beta } \right)} \right)\right| = \frac{\pi}{2}
\]
for any saddle point $t_\pm ^{\left( k \right)}$, we infer that \eqref{eq5} is valid as long as $-\frac{\pi}{2} < \theta < \frac{\pi}{2}$ with a contour $\Gamma\left(\theta\right)$ specified above. When $\theta = -\frac{\pi}{2}$, the path $\mathscr{P}\left(\theta\right)$ connects to the saddle point $t_ + ^{\left( 0 \right)} = i\beta$. Similarly, when $\theta = \frac{\pi}{2}$, the path $\mathscr{P}\left(\theta\right)$ connects to the saddle point $t_ - ^{\left( 0 \right)} = -i\beta$. These are the adjacent saddles. The set
\[
\Delta  = \left\{ {u \in \mathscr{P}\left( \theta  \right) : - \frac{\pi}{2} < \theta  < \frac{\pi}{2}} \right\}
\]
forms a domain in the complex plane whose boundary contains portions of steepest descent paths through the adjacent saddles (see Figure \ref{fig2}). These paths are $\mathscr{C}_ + ^{\left( 0 \right)} \left( { \frac{\pi }{2}} \right)$ and $\mathscr{C}_ - ^{\left( 0 \right)} \left( { - \frac{\pi }{2}} \right)$, and they are called the adjacent contours to the endpoint $0$. The function under the contour integral in \eqref{eq5} is an analytic function of $u$ in the domain $\Delta$, therefore we can deform $\Gamma$ over the adjacent contours. We thus find that for $-\frac{\pi}{2} < \theta < \frac{\pi}{2}$ and $N \geq 0$, \eqref{eq5} may be written
\begin{gather}\label{eq6}
\begin{split}
R_N \left( {\nu ,\beta } \right) = \; & \frac{1}{{\pi \nu ^{2N + 1} }}\int_0^{ + \infty } {s^{2N} e^{ - s} \frac{1}{{2\pi i}}\int_{\mathscr{C}_ + ^{\left( 0 \right)} \left( {\frac{\pi }{2}} \right)} {\frac{{f^{ - 2N - 1} \left( {u,\beta } \right)}}{{1 - \left( {s/\nu } \right)^2 f^{ - 2} \left( {u,\beta } \right)}}du} ds} \\ & + \frac{1}{{\pi \nu ^{2N + 1} }}\int_0^{ + \infty } {s^{2N} e^{ - s} \frac{1}{{2\pi i}}\int_{\mathscr{C}_ - ^{\left( 0 \right)} \left( { - \frac{\pi }{2}} \right)} {\frac{{f^{ - 2N - 1} \left( {u,\beta } \right)}}{{1 - \left( {s/\nu } \right)^2 f^{ - 2} \left( {u,\beta } \right)}}du} ds} .
\end{split}
\end{gather}
Now we make the changes of variable
\[
s = t\frac{{\left| {f\left( {i\beta ,\beta } \right)-f\left(0,\beta\right)} \right|}}{{f\left( {i\beta ,\beta } \right)-f\left(0,\beta\right)}}f\left( {u,\beta } \right) =  - itf\left( {u,\beta } \right)
\]
in the first, and
\[
s = t\frac{{\left| {f\left( { - i\beta ,\beta } \right)-f\left(0,\beta\right)} \right|}}{{f\left( { - i\beta ,\beta } \right)-f\left(0,\beta\right)}}f\left( {u,\beta } \right) = itf\left( {u,\beta } \right)
\]
in the second double integral. Clearly, by the definition of the adjacent contours, $t$ is positive. The quantities $f\left( {i\beta ,\beta } \right) -f\left(0,\beta\right) = i\left( {\tan \beta  - \beta  } \right)$ and $f\left( { - i\beta ,\beta } \right)-f\left(0,\beta\right) = -i\left( {\tan \beta  - \beta } \right)$ were essentially called the ``singulants" by Dingle \cite[p. 147]{Dingle}. With these changes of variable, the representation \eqref{eq6} for $R_N \left( {\nu ,\beta } \right)$ becomes
\begin{equation}\label{eq7}
R_N \left( {\nu ,\beta } \right) = \frac{{\left( { - 1} \right)^N }}{{\pi \nu ^{2N + 1} }}\int_0^{ + \infty } {\frac{{t^{2N} }}{{1 + \left( {t/\nu } \right)^2 }}\left( {\frac{1}{{2\pi }}\int_{\mathscr{C}_ - ^{\left( 0 \right)} \left( { - \frac{\pi }{2}} \right)} {e^{ - itf\left( {u,\beta } \right)} du}  - \frac{1}{{2\pi }}\int_{\mathscr{C}_ + ^{\left( 0 \right)} \left( {\frac{\pi }{2}} \right)} {e^{itf\left( {u,\beta } \right)} du} } \right)dt} ,
\end{equation}
for $-\frac{\pi}{2} < \theta < \frac{\pi}{2}$ and $N \geq 0$. Finally, the contour integrals can themselves be represented in terms of the Hankel functions since
\[
\frac{1}{{2\pi }}\int_{\mathscr{C}_ - ^{\left( 0 \right)} \left( { - \frac{\pi }{2}} \right)} {e^{ - itf\left( {u,\beta } \right)} du}  =  - \frac{i}{2}H_{ - it}^{\left( 2 \right)} \left( { - it\sec \beta } \right) = \frac{i}{2}H_{it}^{\left( 1 \right)} \left( {it\sec \beta } \right),
\]
and
\[
-\frac{1}{{2\pi }}\int_{\mathscr{C}_ + ^{\left( 0 \right)} \left( {\frac{\pi }{2}} \right)} {e^{itf\left( {u,\beta } \right)} du}  =  \frac{i}{2}H_{it}^{\left( 1 \right)} \left( {it\sec \beta } \right) .
\]
Substituting these into \eqref{eq7} gives \eqref{eq8}. To prove the second representation in \eqref{eq9}, we apply \eqref{eq8} for the right-hand side of
\[
a_n \left( { - \sec \beta } \right) = \pi \frac{\nu ^{2n + 1}}{\left( {2n} \right)!}\left( {R_{n + 1} \left( {\nu ,\beta } \right) - R_n \left( {\nu ,\beta } \right)} \right) .
\]

\subsection{Case (ii): $x=1$} We assume that $\theta = 0$ and later we shall use an analytic continuation argument to extend the results to complex $\nu$. Let $f\left( t \right) = \sinh t-t$. If
\begin{equation}\label{eq11}
\tau  = f\left( t \right),
\end{equation}
then $\tau$ is real on the curve $\mathscr{P}\left(0\right)$, and, as $t$ travels along this curve from $0$ to $+\infty $, $\tau$ increases from $0$ to $+\infty$. Therefore, corresponding to each positive value of $\tau$, there is a value of $t$, say $t\left(\tau\right)$, satisfying \eqref{eq11} with $t\left(\tau\right)>0$. In terms of $\tau$, we have
\[
\mathbf{A}_{ - \nu } \left( \nu  \right) = \frac{1}{\pi }\int_0^{ + \infty } {e^{ - \nu \tau } \frac{{dt\left( \tau  \right)}}{{d\tau }}d\tau }  = \frac{1}{\pi }\int_0^{ + \infty } {e^{ - \nu \tau } \frac{1}{{\cosh t\left( \tau  \right) - 1}}d\tau } .
\]
As in the first case, we express the function involving $t\left(\tau\right)$ as a contour integral using the residue theorem, to obtain
\[
\mathbf{A}_{ - \nu } \left( \nu  \right) = \frac{1}{{3\pi }}\int_0^{ + \infty } {\tau ^{ - \frac{2}{3}} e^{ - \nu \tau } \frac{1}{{2\pi i}}\oint_\Gamma  {\frac{{f^{ - \frac{1}{3}} \left( u \right)}}{{1 - \tau ^{\frac{2}{3}} f^{ - \frac{2}{3}} \left( u \right)}}du} d\tau } 
\]
where the contour $\Gamma$ encircles the path $\mathscr{P}\left(0\right)$ in the positive direction and does not enclose any of the saddle points $t ^{\left( k \right)} \neq t ^{\left( 0 \right)}$ (cf. Figure \ref{fig1}). The cube root is defined so that $f^{\frac{1}{3}}\left( t \right)$ is positive on the path $\mathscr{P}\left(0\right)$. Next we apply the expression \eqref{eq2} to expand the function under the contour integral in powers of $\tau^{\frac{2}{3}} f^{-\frac{2}{3}} \left( u \right)$. The result is
\[
\mathbf{A}_{ - \nu } \left( \nu  \right) = \frac{1}{{3\pi }}\sum\limits_{n = 0}^{N - 1} {\int_0^{ + \infty } {\tau ^{\frac{{2n - 2}}{3}} e^{ - \nu \tau } \frac{1}{{2\pi i}}\oint_\Gamma  {\frac{{du}}{{f^{\frac{{2n + 1}}{3}} \left( u \right)}}} d\tau } }  + R_N \left( \nu  \right)
\]
where
\begin{equation}\label{eq12}
R_N \left( \nu  \right) = \frac{1}{{3\pi }}\int_0^{ + \infty } {\tau ^{\frac{{2N - 2}}{3}} e^{ - \nu \tau } \frac{1}{{2\pi i}}\oint_\Gamma  {\frac{{f^{ - \frac{{2N + 1}}{3}} \left( u \right)}}{{1 - \tau ^{\frac{2}{3}} f^{ - \frac{2}{3}} \left( u \right)}}du} d\tau } .
\end{equation}
The path $\Gamma$ in the sum can be shrunk into a small circle around $t^{\left( 0 \right)} = 0$, and we arrive at
\[
\mathbf{A}_{ - \nu } \left( \nu  \right) = \frac{1}{{3\pi }}\sum\limits_{n = 0}^{N - 1} {d_{2n} \frac{{\Gamma \left( {\frac{{2n + 1}}{3}} \right)}}{{\nu ^{\frac{{2n + 1}}{3}} }}}  + R_N \left( \nu  \right)
\]
where
\[
d_{2n}  = \frac{1}{{2\pi i}}\oint_{\left( {0^ +  } \right)} {\frac{{du}}{{f^{\frac{{2n + 1}}{3}} \left( u \right)}}}  = \frac{1}{{\left( {2n} \right)!}}\left[ {\frac{{d^{2n} }}{{dt^{2n} }}\left( {\frac{{t^3 }}{{\sinh t - t}}} \right)^{\frac{{2n + 1}}{3}} } \right]_{t = 0} .
\]

\begin{figure}[!t]
\def\svgwidth{0.6\columnwidth}
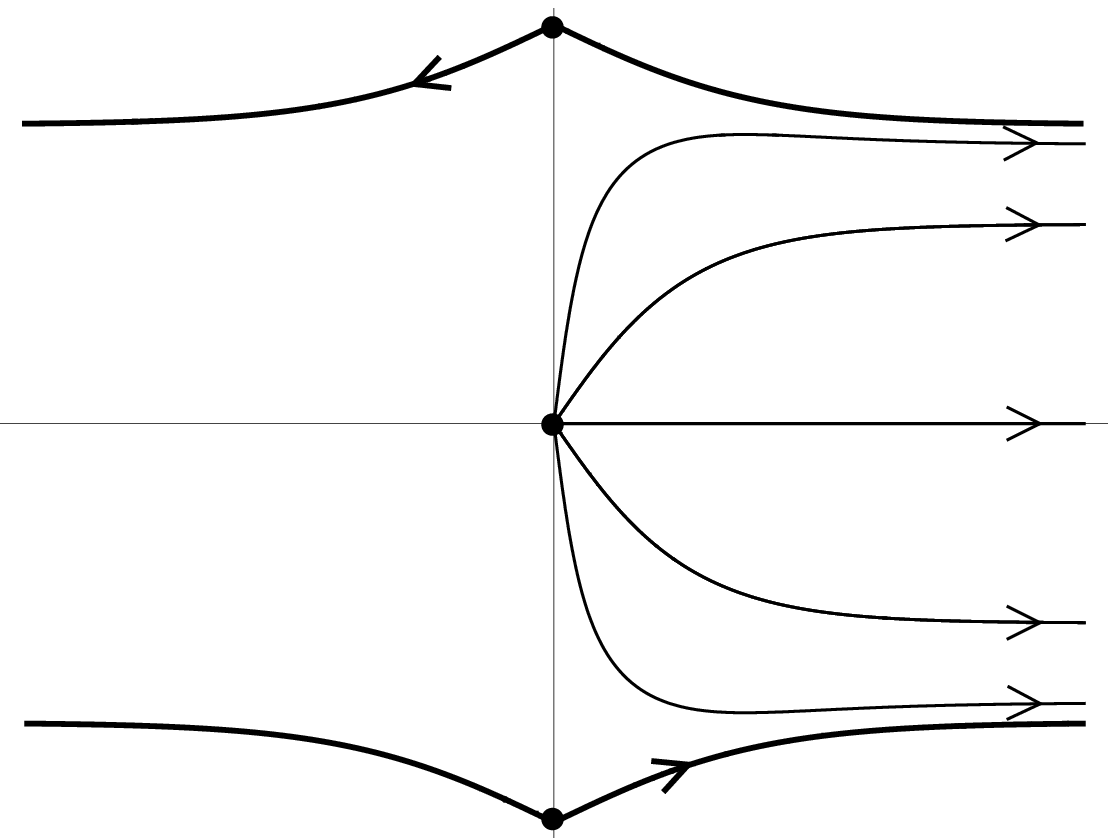
\caption{The path $\mathscr{P}\left(\theta\right)$ emanating from the saddle point $t^{\left(0\right)}$ when (i) $\theta=0$, (ii) $\theta=-\pi$, (iii) $\theta=-\frac{7\pi}{5}$, (iv) $\theta=\pi$, (v) $\theta=\frac{7\pi}{5}$. The paths $\mathscr{L}^{\left( 1 \right)} \left( { -\frac{\pi }{2}} \right)$ and $\mathscr{P} ^{\left( -1 \right)} \left( { \frac{\pi }{2}} \right)$ are the adjacent contours to $t^{\left(0\right)}$. The domain $\Delta$ comprises all points between these two paths in the right-half plane.}
\label{fig3}
\end{figure}

Applying the change of variable $\nu \tau = s$ in \eqref{eq12} gives
\begin{equation}\label{eq13}
R_N \left( \nu  \right) = \frac{1}{{3\pi \nu ^{\frac{{2N + 1}}{3}} }}\int_0^{ + \infty } {s^{\frac{{2N - 2}}{3}} e^{ - s} \frac{1}{{2\pi i}}\oint_\Gamma  {\frac{{f^{ - \frac{{2N + 1}}{3}} \left( u \right)}}{{1 - \left( {s/\nu } \right)^{\frac{2}{3}} f^{ - \frac{2}{3}} \left( u \right)}}du} ds} .
\end{equation}
As in the first case, we need to locate the adjacent saddle points. When $\theta = -\frac{3\pi}{2}$, the path $\mathscr{P}\left(\theta\right)$ connects to the saddle point $t ^{\left( 1 \right)} = 2\pi i$. Similarly, when $\theta = \frac{3\pi}{2}$, the path $\mathscr{P}\left(\theta\right)$ connects to the saddle point $t ^{\left( -1 \right)} = -2\pi i$. Therefore, the adjacent saddles are $t^{\left(\pm 1\right)}$. The set
\[
\Delta  = \left\{ {u \in \mathscr{P}\left( \theta  \right) : - \frac{3\pi}{2} < \theta  < \frac{3\pi}{2}} \right\}
\]
forms a domain in the complex plane whose boundary contains portions of steepest descent paths through the adjacent saddles (see Figure \ref{fig3}). These paths are $\mathscr{L} ^{\left( 1 \right)} \left( { -\frac{\pi }{2}} \right)$ and $\mathscr{P}^{\left( -1 \right)} \left( { \frac{\pi }{2}} \right)$, the adjacent contours to the saddle point $t ^{\left( 0 \right)}$ (these paths are defined in \cite{Nemes}). The function under the contour integral in \eqref{eq13} is an analytic function of $u$ in the domain $\Delta$, therefore we can deform $\Gamma$ over the adjacent contours. We thus find that for $-\frac{3\pi}{2} < \theta < \frac{3\pi}{2}$ and $N \geq 0$, \eqref{eq13} may be written
\begin{gather}\label{eq14}
\begin{split}
R_N \left( \nu  \right) = \; & \frac{1}{{3\pi \nu ^{\frac{{2N + 1}}{3}} }}\int_0^{ + \infty } {s^{\frac{{2N - 2}}{3}} e^{ - s} \frac{1}{{2\pi i}}\int_{\mathscr{L}^{\left( 1 \right)} \left( { - \frac{\pi }{2}} \right)} {\frac{{f^{ - \frac{{2N + 1}}{3}} \left( u \right)}}{{1 - \left( {s/\nu } \right)^{\frac{2}{3}} f^{ - \frac{2}{3}} \left( u \right)}}du} ds} \\ & + \frac{1}{{3\pi \nu ^{\frac{{2N + 1}}{3}} }}\int_0^{ + \infty } {s^{\frac{{2N - 2}}{3}} e^{ - s} \frac{1}{{2\pi i}}\int_{\mathscr{P}^{\left( { - 1} \right)} \left( {\frac{\pi }{2}} \right)} {\frac{{f^{ - \frac{{2N + 1}}{3}} \left( u \right)}}{{1 - \left( {s/\nu } \right)^{\frac{2}{3}} f^{ - \frac{2}{3}} \left( u \right)}}du} ds} .
\end{split}
\end{gather}
Now we perform the changes of variable
\[
s = t\frac{{\left| {f\left( {2\pi i} \right) - f\left( 0 \right)} \right|}}{{f\left( 2\pi i \right) - f\left( 0 \right)}} f\left( u \right) = itf\left( u \right)
\]
in the first, and
\[
s = t\frac{{\left| {f\left( { - 2\pi i} \right) - f\left( 0 \right)} \right|}}{{f\left( { - 2\pi i} \right) - f\left( 0 \right)}} f\left( u \right) =  - itf\left( u \right)
\]
in the second double integral. In this case, Dingle's singulants are $f\left( { \pm 2\pi i} \right) - f\left( 0 \right) = \mp 2\pi i$. When using these changes of variable, we should take $i^{\frac{2}{3}}  = -1$ in the first, and $\left( { - i} \right)^{\frac{2}{3}}  = -1$ in the second double integral. With these changes of variable, the representation \eqref{eq14} for $R_N \left( \nu  \right)$ becomes
\begin{equation}\label{eq15}
R_N \left( \nu  \right) = \frac{{\left( { - 1} \right)^N }}{{3\pi \nu ^{\frac{{2N + 1}}{3}} }}\int_0^{ + \infty } {\frac{{t^{\frac{{2N - 2}}{3}} }}{{1 + \left( {t/\nu } \right)^{\frac{2}{3}} }}\left( {\frac{1}{{2\pi }}\int_{\mathscr{P}^{\left( { - 1} \right)} \left( {\frac{\pi }{2}} \right)} {e^{itf\left( u \right)} du}  - \frac{1}{{2\pi }}\int_{\mathscr{L}^{\left( 1 \right)} \left( { - \frac{\pi }{2}} \right)} {e^{ - itf\left( u \right)} du} } \right)dt} ,
\end{equation}
for $-\frac{3\pi}{2} < \theta < \frac{3\pi}{2}$ and $N \geq 0$. Finally, the contour integrals can themselves be represented in terms of the Hankel functions since
\[
\frac{1}{{2\pi }}\int_{\mathscr{P}^{\left( { - 1} \right)} \left( {\frac{\pi }{2}} \right)} {e^{itf\left( u \right)} du}  = \frac{{e^{ - 2\pi t} }}{2}iH_{it}^{\left( 1 \right)} \left( {it} \right),
\]
and
\[
 - \frac{1}{{2\pi }}\int_{\mathscr{L}^{\left( 1 \right)} \left( { - \frac{\pi }{2}} \right)} {e^{ - itf\left( u \right)} du} = - \frac{{e^{ - 2\pi t} }}{2}iH_{ - it}^{\left( 2 \right)} \left( { - it} \right) =  \frac{{e^{ - 2\pi t} }}{2}iH_{it}^{\left( 1 \right)} \left( {it} \right) .
\]
Substituting these into \eqref{eq15} gives \eqref{eq16}. To prove the second representation in \eqref{eq17}, we apply \eqref{eq16} for the right-hand side of
\[
d_{2n}  = 3\pi \frac{{\nu ^{\frac{{2n + 1}}{3}} }}{{\Gamma \left( {\frac{{2n + 1}}{3}} \right)}}\left( {R_n \left( \nu  \right) - R_{n + 1} \left( \nu  \right)} \right).
\]

\section{Error bounds}\label{section3}

In this section we derive explicit and realistic error bounds for the large order asymptotic series of the Anger--Weber function. The proofs are based on the resurgence formulas given in Theorems \ref{thm1} and \ref{thm2}.

We comment on the relation between Meijer's work on the asymptotic expansion of $\mathbf{A}_{ - \nu } \left( {\nu \sec \beta } \right)$ \cite{Meijer} and ours. Some of the estimates in \cite{Meijer} coincide with ours and are valid in wider sectors of the complex $\nu$-plane. However, it should be noted that those bounds become less effective outside the sectors of validity of the representation \eqref{eq8} due to the Stokes phenomenon. For those sectors we recommend the use of the continuation formulas given in Section \ref{section1}.

To estimate the remainder terms, we shall use the elementary result that
\begin{equation}\label{eq20}
\frac{1}{{\left| {1 - re^{i\varphi } } \right|}} \le \begin{cases} \left|\csc \varphi \right| & \; \text{ if } \; 0 < \left|\varphi \text{ mod } 2\pi\right| <\frac{\pi}{2} \\ 1 & \; \text{ if } \; \frac{\pi}{2} \leq \left|\varphi \text{ mod } 2\pi\right| \leq \pi \end{cases}
\end{equation}
holds for any $r>0$. We will also need the fact that
\begin{equation}\label{eq21}
iH_{it}^{\left(1\right)} \left( {itx} \right) \ge 0
\end{equation}
for any $t>0$ and $x\geq 1$ (see \cite{Nemes}).

\subsection{Case (i): $x>1$} As usual, let $0<\beta<\frac{\pi}{2}$ be defined by $\sec \beta = x$. We observe that from \eqref{eq9} and \eqref{eq21} it follows that
\[
\left| {a_n \left( { - \sec \beta } \right)} \right| = \frac{1}{{\left( {2n} \right)!}}\int_0^{ + \infty } {t^{2n} iH_{it}^{\left( 1 \right)} \left( {it\sec \beta } \right)dt} .
\]
Using this formula, together with the representation \eqref{eq8} and the estimate \eqref{eq20}, we obtain the error bound
\begin{equation}\label{eq23}
\left| {R_N \left( {\nu ,\beta } \right)} \right| \le \frac{1}{\pi }\frac{{\left( {2N} \right)!\left| {a_N \left( { - \sec \beta } \right)} \right|}}{{\left| \nu  \right|^{2N + 1} }} \begin{cases} \left|\csc\left(2\theta\right)\right| & \; \text{ if } \; \frac{\pi}{4} < \left|\theta\right| <\frac{\pi}{2} \\ 1 & \; \text{ if } \; \left|\theta\right| \leq \frac{\pi}{4}. \end{cases}
\end{equation}
Here and throughout, $\theta = \arg \nu$. When $\nu$ is real and positive, we can obtain more precise estimates. Indeed, as $0 < \frac{1}{1 + \left( {t/\nu } \right)^2} < 1 $ for $t,\nu>0$, from \eqref{eq8} and \eqref{eq9} we find
\[
R_N \left( {\nu ,\beta } \right) =  - \frac{1}{\pi }\frac{{\left( {2N} \right)!a_N \left( { - \sec \beta } \right)}}{{\nu ^{2N + 1} }}\Theta ,
\]
where $0 < \Theta < 1$ is an appropriate number depending on $\nu,\beta$ and $N$. In particular, when $N=0$, we have
\[
0 < \mathbf{A}_{ - \nu } \left( {\nu \sec \beta } \right) < \frac{1}{{\pi \nu \left( {\sec \beta  - 1} \right)}} \; \text{ for } \; \nu >0.
\]
Therefore, the leading order asymptotic approximation for $\mathbf{A}_{ - \nu } \left( \nu \sec \beta \right)$ is always in error by excess, for all positive values of $\nu$ (cf. \cite[p. 298, formula 11.11.14]{NIST}).

The error bound \eqref{eq23} becomes singular as $\theta \to \pm \frac{\pi}{2}$, and therefore unrealistic near the Stokes lines. A better bound for $R_N \left( {\nu ,\beta } \right)$ near these lines can be derived as follows. Let $0 < \varphi  < \frac{\pi }{2}$ be an acute angle that may depend on $N$. Suppose that $\frac{\pi}{4} +\varphi < \theta \le \frac{\pi}{2}$. An analytic continuation of the representation \eqref{eq18} to this sector can be found by rotating the path of integration in \eqref{eq8} by $\varphi$:
\[
R_N \left( {\nu ,\beta } \right) = \frac{{\left( { - 1} \right)^N }}{{\pi \nu ^{2N + 1} }}\int_0^{ + \infty e^{i\varphi } } {\frac{{t^{2N} }}{{1 + \left( {t/\nu } \right)^2 }}iH_{it}^{\left( 1 \right)} \left( {it\sec \beta } \right)dt} .
\]
Substituting $t = \frac{se^{i\varphi }}{\cos \varphi}$ and applying the estimation \eqref{eq20}, we obtain
\[
\left|R_N \left( {\nu ,\beta } \right)\right| \le \frac{{\csc \left( {2\left( {\theta  - \varphi } \right)} \right)}}{{\pi \cos ^{2N + 1} \varphi \left| \nu  \right|^{2N + 1} }}\int_0^{ + \infty } {s^{2N} \left| {H_{\frac{{ise^{i\varphi } }}{{\cos \varphi }}}^{\left( 1 \right)} \left( {\frac{{ise^{i\varphi } }}{{\cos \varphi }}\sec \beta } \right)} \right|ds} .
\]
In \cite{Nemes}, it was shown that
\begin{equation}\label{eq27}
\left| {H_{\frac{{ise^{i\varphi } }}{{\cos \varphi }}}^{\left( 1 \right)} \left( {\frac{{ise^{i\varphi } }}{{\cos \varphi }}\sec \beta } \right)} \right| \le \frac{1}{{\sqrt {\cos \varphi } }}\left| {H_{is}^{\left( 1 \right)} \left( {is\sec \beta } \right)} \right| = \frac{1}{{\sqrt {\cos \varphi } }}iH_{is}^{\left( 1 \right)} \left( {is\sec \beta } \right)
\end{equation}
for any $s>0$ and $0 < \varphi  < \frac{\pi }{2}$. It follows that
\begin{gather}\label{eq24}
\begin{split}
\left|R_N \left( {\nu ,\beta } \right)\right| & \le \frac{{\csc \left( {2\left( {\theta  - \varphi } \right)} \right)}}{{\pi \cos ^{2N + \frac{3}{2}} \varphi \left| \nu  \right|^{2N + 1} }}\int_0^{ + \infty } {s^{2N} iH_{is}^{\left( 1 \right)} \left( {is\sec \beta } \right)ds} \\
& = \frac{{\csc \left( {2\left( {\theta  - \varphi } \right)} \right)}}{{\cos ^{2N + \frac{3}{2}} \varphi }}\frac{1}{\pi }\frac{{\left( {2N} \right)!\left| {a_N \left( { - \sec \beta } \right)} \right|}}{{\left| \nu  \right|^{2N + 1} }} .
\end{split}
\end{gather}
The angle $\varphi  = \arctan \left( {\left( {\frac{{4N + 5}}{2}} \right)^{ - \frac{1}{2}} } \right)$ minimizes the function $\csc \left( {2\left( {\frac{\pi }{2} - \varphi } \right)} \right)\cos ^{ - 2N - \frac{3}{2}} \varphi$, and
\begin{align*}
\frac{{\csc \left( {2\left( {\theta  - \arctan \left( {\left( {\frac{{4N + 5}}{2}} \right)^{ - \frac{1}{2}} } \right)} \right)} \right)}}{{\cos ^{2N + \frac{3}{2}} \left( {\arctan \left( {\left( {\frac{{4N + 5}}{2}} \right)^{ - \frac{1}{2}} } \right)} \right)}} & \le \frac{{\csc \left( {2\left( {\frac{\pi }{2} - \arctan \left( {\left( {\frac{{4N + 5}}{2}} \right)^{ - \frac{1}{2}} } \right)} \right)} \right)}}{{\cos ^{2N + \frac{3}{2}} \left( {\arctan \left( {\left( {\frac{{4N + 5}}{2}} \right)^{ - \frac{1}{2}} } \right)} \right)}} \\ &= \frac{1}{{\sqrt 2 }}\left( 1+\frac{2}{4N + 5} \right)^{N + \frac{7}{4}} \sqrt {N + \frac{5}{4}}  \le \sqrt {\frac{e}{2}\left( {N + \frac{3}{2}} \right)}
\end{align*}
for all $\frac{\pi }{4} + \varphi  = \frac{\pi }{4} + \arctan \left( {\left( {\frac{{4N + 5}}{2}} \right)^{ - \frac{1}{2}} } \right) < \theta  \le \frac{\pi }{2}$ with $N \geq 0$. Applying this in \eqref{eq24} yields the upper bound
\begin{equation}\label{eq25}
\left| {R_N \left( {\nu ,\beta } \right)} \right| \le \sqrt {\frac{e}{2}\left( {N + \frac{3}{2}} \right)} \frac{1}{\pi }\frac{{\left( {2N} \right)!\left| {a_N \left( { - \sec \beta } \right)} \right|}}{{\left| \nu  \right|^{2N + 1} }},
\end{equation}
which is valid for $\frac{\pi }{4} + \varphi  = \frac{\pi }{4} + \arctan \left( {\left( {\frac{{4N + 5}}{2}} \right)^{ - \frac{1}{2}} } \right) < \theta  \le \frac{\pi }{2}$ with $N \geq 0$. Since $\left| {R_N \left( {\bar \nu ,\beta } \right)} \right| = \left| {\overline {R_N \left( {\nu ,\beta } \right)} } \right| = \left| {R_N \left( {\nu ,\beta } \right)} \right|$, this bound also holds when $-\frac{\pi}{2} \leq \theta < -\frac{\pi}{4} - \arctan \left( {\left( {\frac{{4N + 5}}{2}} \right)^{ - \frac{1}{2}} } \right)$. In the ranges $\frac{\pi }{4} < \left| \theta  \right| \leq \frac{\pi }{4} + \arctan \left( {\frac{{\sqrt 2 }}{3}} \right)$ it holds that $\left| {\csc \left( {2\theta } \right)} \right| \le \sqrt {\frac{e}{2}\left( {1 + \frac{3}{2}} \right)}$, whence the estimate \eqref{eq25} is valid in the wider sectors $\frac{\pi }{4} < \left| \theta  \right| \le \frac{\pi }{2}$ as long as $N\geq 1$.

\subsection{Case (ii): $x=1$} We note that from \eqref{eq17} and \eqref{eq21} it follows that
\[
\left| {d_{2n} } \right| = \frac{1}{{\Gamma \left( {\frac{{2n + 1}}{3}} \right)}}\int_0^{ + \infty } {t^{\frac{{2n - 2}}{3}} e^{ - 2\pi t} iH_{it}^{\left( 1 \right)} \left( {it} \right)dt} .
\]
Applying this formula together with the representation \eqref{eq16} and the inequality \eqref{eq20} yields the error bound
\[
\left| {R_N \left( \nu  \right)} \right| \le \frac{1}{{3\pi }}\left| {d_{2N} } \right|\frac{{\Gamma \left( {\frac{{2N + 1}}{3}} \right)}}{{\left| \nu  \right|^{\frac{{2N + 1}}{3}} }} \begin{cases} \left| {\csc \left( {\frac{2}{3}\theta } \right)} \right| & \; \text{ if } \; \frac{3\pi}{4} < \left| \theta  \right| < \frac{3\pi}{2} \\ 1 & \; \text{ if } \; \left| \theta  \right| \le \frac{3\pi}{4}. \end{cases}
\]
Again, when $\nu$ is real and positive, we can deduce better estimates. Indeed, as $0 < \frac{1}{1 + \left( {t/\nu } \right)^{\frac{2}{3} }} < 1 $ for $t,\nu>0$, from \eqref{eq16} and \eqref{eq17} we find
\[
R_N \left( \nu  \right) = \frac{1}{{3\pi }}d_{2N} \frac{{\Gamma \left( {\frac{{2N + 1}}{3}} \right)}}{{\nu ^{\frac{{2N + 1}}{3}} }}\Xi ,
\]
where $0 < \Xi < 1$ is a suitable number depending on $\nu$ and $N$. In particular, when $N=0$, we have
\[
0 < \mathbf{A}_{ - \nu } \left( \nu \right) < \frac{1}{{3\pi }}d_0 \frac{{\Gamma \left( {\frac{1}{3}} \right)}}{{\nu ^{\frac{1}{3}} }} = \frac{{2^{\frac{4}{3}} }}{{3^{\frac{7}{6}} \Gamma \left( {\frac{2}{3}} \right)\nu ^{\frac{1}{3}} }} \; \text{ for } \; \nu >0.
\]
Hence, the leading order asymptotic approximation for $\mathbf{A}_{ - \nu } \left( \nu \right)$ is always in error by excess, for all positive values of $\nu$ (cf. \cite[p. 298, formula 11.11.16]{NIST}).

Our bound for $R_N \left( \nu \right)$ is unrealistic near the Stokes lines $\theta = \pm\frac{3\pi }{2}$ due to the presence of the factor $\csc \left( {\frac{2}{3}\theta } \right)$. We shall derive better bounds for $R_N \left( \nu \right)$ near these lines using the method we applied in the previous case. Let $0 < \varphi  < \frac{\pi }{2}$ be an acute angle that may depend on $N$ and suppose that $\frac{3\pi}{4} + \varphi  < \theta  \le \frac{3\pi}{2}$. We rotate the path of integration in \eqref{eq16} by $\varphi$, and apply the inequality \eqref{eq20} to obtain
\[
\left| {R_N \left( \nu  \right)} \right| \le \frac{{\csc \left( {\frac{2}{3}\left( {\theta  - \varphi } \right)} \right)}}{{3\pi \cos ^{\frac{{2N + 1}}{3}} \varphi \left| \nu  \right|^{\frac{{2N + 1}}{3}} }}\int_0^{ + \infty } {s^{\frac{{2N - 2}}{3}} e^{ - 2\pi s} \left| {H_{\frac{{ise^{i\varphi } }}{{\cos \varphi }}}^{\left( 1 \right)} \left( {\frac{{ise^{i\varphi } }}{{\cos \varphi }}} \right)} \right|ds} 
\]
for $\frac{3\pi}{4} + \varphi  < \theta  \le \frac{3\pi}{2}$ and $N\geq 0$. Using a continuity argument for the inequality \eqref{eq27}, yields
\[
\left| {H_{\frac{{it}}{{e^{i\varphi } \cos \varphi }}}^{\left( 1 \right)} \left( {\frac{{it}}{{e^{i\varphi } \cos \varphi }}} \right)} \right| \le \frac{1}{{\sqrt {\cos \varphi } }}iH_{it}^{\left( 1 \right)} \left( {it} \right) \le \frac{1}{{\cos ^{\frac{2}{3}} \varphi }}iH_{it}^{\left( 1 \right)} \left( {it} \right) 
\]
for $s>0$ and $0 < \varphi  < \frac{\pi }{2}$. It follows that
\begin{equation}\label{eq28}
\left| {R_N \left( \nu  \right)} \right|  \le \frac{{\csc \left( {\frac{2}{3}\left( {\theta  - \varphi } \right)} \right)}}{{3\pi \cos ^{\frac{{2N + 1}}{3}} \varphi \left| \nu  \right|^{\frac{{2N + 1}}{3}} }}\int_0^{ + \infty } {s^{\frac{{2N - 2}}{3}} e^{ - 2\pi s} i H_{is}^{\left( 1 \right)} \left( is \right)ds} = \frac{\csc \left( {\frac{2}{3}\left( {\theta  - \varphi } \right)} \right)}{ \cos ^{\frac{{2N + 3}}{3}} \varphi} \frac{1}{{3\pi }}\left| {d_{2N} } \right|\frac{{\Gamma \left( {\frac{{2N + 1}}{3}} \right)}}{{\left| \nu  \right|^{\frac{{2N + 1}}{3}} }}.
\end{equation}
There is no simple way to minimize $\csc \left( {\frac{2}{3}\left( {\frac{3\pi}{2}  - \varphi } \right)} \right)\cos ^{ - \frac{{2N + 3}}{3}} \varphi$ in $\varphi$. Nevertheless, an approximate minimizer is given by $\varphi  = \arctan \left( {\left( {\frac{{2N + 2}}{3}} \right)^{ - \frac{1}{2}} } \right)$. It is elementary to show that
\[
\frac{{\csc \left( {\frac{2}{3}\left( {\theta  - \arctan \left( {\left( {\frac{{2N + 2}}{3}} \right)^{ - \frac{1}{2}} } \right) } \right)} \right)}}{{\cos ^{\frac{{2N + 3}}{3}} \left(\arctan \left( {\left( {\frac{{2N + 2}}{3}} \right)^{ - \frac{1}{2}} } \right)\right) }} \le \frac{{\csc \left( {\frac{2}{3}\left( {\frac{{3\pi }}{2} - \arctan \left( {\left( {\frac{{2N + 2}}{3}} \right)^{ - \frac{1}{2}} } \right) } \right)} \right)}}{{\cos ^{\frac{{2N + 3}}{3}} \left(\arctan \left( {\left( {\frac{{2N + 2}}{3}} \right)^{ - \frac{1}{2}} } \right)\right) }} \le \sqrt {\frac{3e}{2}\left( {N + 2} \right)} 
\]
for $\frac{3\pi}{4} + \varphi = \frac{3\pi}{4} + \arctan \left( {\left( {\frac{{2N + 2}}{3}} \right)^{ - \frac{1}{2}} } \right)  < \theta  \le \frac{3\pi}{2}$ and $N\geq 0$. Employing this estimate in \eqref{eq28} gives the upper bound
\begin{equation}\label{eq29}
\left| {R_N \left( \nu  \right)} \right| \le \sqrt {\frac{3e}{2}\left( {N + 2} \right)} \frac{1}{{3\pi }}\left| {d_{2N} } \right|\frac{{\Gamma \left( {\frac{{2N + 1}}{3}} \right)}}{{\left| \nu  \right|^{\frac{{2N + 1}}{3}} }},
\end{equation}
valid when $\frac{3\pi}{4} + \varphi = \frac{3\pi}{4} + \arctan \left( {\left( {\frac{{2N + 2}}{3}} \right)^{ - \frac{1}{2}} } \right)  < \theta  \le \frac{3\pi}{2}$ and $N\geq 0$. A similar argument shows that this bound also holds in the sector $-\frac{3\pi}{2} \leq \theta < -\frac{3\pi}{4} - \arctan \left( {\left( {\frac{{2N + 2}}{3}} \right)^{ - \frac{1}{2}} } \right)$. In the ranges $\frac{{3\pi }}{4} < \left| \theta  \right| \le \frac{{3\pi }}{4} + \arctan \left( {\left( {\frac{2}{3}} \right)^{ - \frac{1}{2}} } \right)$ it holds that $\left| {\csc \left( {\frac{2}{3}\theta } \right)} \right| \le \sqrt {\frac{3e}{2}\left( {0 + 2} \right)}  = \sqrt {3e}$, therefore, the estimate \eqref{eq29} remains valid in the wider sectors $\frac{3\pi}{4} < \left| \theta  \right| \le \frac{3\pi}{2}$ for any $N\geq 0$.

\section{Asymptotics for the late coefficients}\label{section4}

In this section, we investigate the asymptotic nature of the coefficients $a_n\left(- \sec \beta\right)$ as $n \to +\infty$. The asymptotic behaviour of the coefficients $d_{2n}$ is discussed in the earlier paper \cite{Nemes}. For our purposes, the most appropriate representation of the coefficients $a_n\left(- \sec \beta\right)$ is the second integral formula in \eqref{eq9}. From \eqref{eq40}, it follows that for any $t>0$ and $0 < \beta <\frac{\pi}{2}$, it holds that
\begin{equation}\label{eq31}
iH_{it}^{\left( 1 \right)} \left( {it\sec \beta } \right) = \frac{{e^{ - t\left( {\tan \beta  - \beta } \right)} }}{{\left( {\frac{1}{2}t\pi \tan \beta } \right)^{\frac{1}{2}} }}\left( {\sum\limits_{m = 0}^{M - 1} {\frac{{i^m U_m \left( {i\cot \beta } \right)}}{{t^m }}}  + R_M^{\left( H \right)} \left( {it,\beta } \right)} \right).
\end{equation}
In \cite{Nemes}, it was proved that the remainder $R_M^{\left( H \right)} \left( {it,\beta } \right)$ satisfies
\begin{equation}\label{eq32}
\left| {R_M^{\left( H \right)} \left( {it,\beta } \right)} \right| \le \frac{{\left| {U_M \left( {i\cot \beta } \right)} \right|}}{{t^M }} .
\end{equation}
Substituting the formula \eqref{eq31} into \eqref{eq9} gives the expansion
\begin{gather}\label{eq30}
\begin{split}
 - \left( {2n} \right)!a_n \left( { - \sec \beta } \right) = & \left( {\frac{{2\cot \beta }}{{\pi \left( {\tan \beta  - \beta } \right)}}} \right)^{\frac{1}{2}} \frac{{\left( { - 1} \right)^n \Gamma \left( {2n + \frac{1}{2}} \right)}}{{\left( {\tan \beta  - \beta } \right)^{2n} }} \\ & \times \left( {\sum\limits_{m = 0}^{M - 1} {\left( {i\left( {\tan \beta  - \beta } \right)} \right)^m U_m \left( {i\cot \beta } \right)\frac{{\Gamma \left( {2n - m + \frac{1}{2}} \right)}}{{\Gamma \left( {2n + \frac{1}{2}} \right)}}}  + A_M \left( {n,\beta } \right)} \right),
\end{split}
\end{gather}
for any fixed $0 \le M \le 2n$, provided that $n\geq 1$. The remainder term $A_M \left( {n,\beta } \right)$ is given by the integral formula
\[
A_M \left( {n,\beta } \right) = \frac{{\left( {\tan \beta  - \beta } \right)^{2n + \frac{1}{2}} }}{{\Gamma \left( {2n + \frac{1}{2}} \right)}}\int_0^{ + \infty } {t^{2n - \frac{1}{2}} e^{ - t\left( {\tan \beta  - \beta } \right)} R_M^{\left( H \right)} \left( {it,\beta } \right)dt} .
\]
To bound this error term, we apply the estimate \eqref{eq32} to find
\begin{equation}\label{eq33}
\left| {A_M \left( {n,\beta } \right)} \right| \le \left( {\tan \beta  - \beta } \right)^M \left| {U_M \left( {i\cot \beta } \right)} \right|\frac{{\Gamma \left( {2n - M + \frac{1}{2}} \right)}}{{\Gamma \left( {2n + \frac{1}{2}} \right)}} .
\end{equation}
Expansions of type \eqref{eq30} are called inverse factorial series in the literature. Numerically, their character is similar to the character of asymptotic power series, because the consecutive Gamma functions decrease asymptotically by a factor $2n$.

From the asymptotic behaviour of the coefficients $U_m \left( {i\cot \beta } \right)$ (see \cite{Nemes}), we infer that for large $n$, the least value of the bound \eqref{eq33} occurs when $M \approx \frac{4n}{3}$. With this choice of $M$, the error bound is $\mathcal{O}\left( {n^{\frac{1}{2}} 9^{-n} } \right)$. This is the best accuracy we can achieve using the expansion \eqref{eq30}.

By extending the sum in \eqref{eq30} to infinity, we arrive at the formal series
\begin{multline*}
 - \left( {2n} \right)!a_n \left( { - \sec \beta } \right) \approx \left( {\frac{{2\cot \beta }}{{\pi \left( {\tan \beta  - \beta } \right)}}} \right)^{\frac{1}{2}} \frac{{\left( { - 1} \right)^n \Gamma \left( {2n + \frac{1}{2}} \right)}}{{\left( {\tan \beta  - \beta } \right)^{2n} }}\left( {1 - \frac{{\left( {\tan \beta  - \beta } \right)\cot \beta \left( {5\cot ^2 \beta  + 3} \right)}}{{24\left( {2n - \frac{1}{2}} \right)}}  }\right. \\
\left.{ + \frac{{\left( {\tan \beta  - \beta } \right)^2 \cot ^2 \beta \left( {385\cot ^4 \beta  + 462\cot ^2 \beta  + 81} \right)}}{{1152\left( {2n - \frac{1}{2}} \right)\left( {2n - \frac{3}{2}} \right)}} +  \cdots } \right).
\end{multline*}
This is exactly Dingle's expansion for the late coefficients in the asymptotic series of $\mathbf{A}_{-\nu} \left(\nu \sec \beta\right)$ \cite[p. 202]{Dingle}. The mathematically rigorous form of Dingle's series is therefore the formula \eqref{eq30}.

Numerical examples illustrating the efficacy of the expansion \eqref{eq30}, truncated optimally, are given in Table \ref{table1}.

\begin{table*}[!ht]
\begin{center}
\begin{tabular}
[c]{ l r @{\,}c@{\,} l}\hline
 & \\ [-1ex]
 values of $\beta$ and $M$ & $\beta=\frac{\pi}{6}$, $M=33$ & & \\ [1ex]
 exact numerical value of $a_{25}\left(-\sec\beta\right)$ & $0.19289505370609710328176787542524$ & $\times$ & $10^{64}$ \\ [1ex]
 approximation \eqref{eq30} to $a_{25}\left(-\sec\beta\right)$ & $0.19289505370609710328176788499115$ & $\times$ & $10^{64}$  \\ [1ex]
 error & $-0.956591$ & $\times$ & $10^{38}$\\ [1ex] 
  error bound using \eqref{eq33} & $0.1871709$ & $\times$ & $10^{39}$\\ [1ex] \hline
 & \\ [-1ex]
 values of $\beta$ and $M$ &  $\beta=\frac{\pi}{3}$, $M=33$ & & \\ [1ex]
 exact numerical value of $a_{25}\left(-\sec\beta\right)$ & $0.17129537192362280172104021636215$ & $\times$ & $10^8$ \\ [1ex]
 approximation \eqref{eq30} to $a_{25}\left(-\sec\beta\right)$ & $0.17129537192362280172104022485431$ & $\times$ & $10^8$  \\ [1ex]
 error & $-0.849216$ & $\times$ & $10^{-18}$\\ [1ex]
 error bound using \eqref{eq33} & $0.1661627$ & $\times$ & $10^{-17}$\\ [1ex] \hline
 & \\ [-1ex]
 values of $\beta$ and $M$ & $\beta=\frac{6\pi}{13}$, $M=33$ & & \\ [1ex]
 exact numerical value of $a_{25}\left(-\sec\beta\right)$ & $0.39520964363504437817499204430357$ & $\times$ & $10^{-43}$ \\ [1ex]
 approximation \eqref{eq30} to $a_{25}\left(-\sec\beta\right)$ & $0.39520964363504437817499206387318$ & $\times$ & $10^{-43}$ \\ [1ex]
 error & $-0.1956961$ & $\times$ & $10^{-68}$\\ [1ex]
 error bound using \eqref{eq33} & $0.3829199$ & $\times$ & $10^{-68}$\\ [1ex] \hline
 & \\ [-1ex]
 values of $\beta$ and $M$ &  $\beta=\frac{7\pi}{15}$, $M=33$ & & \\ [1ex]
 exact numerical value of $a_{25}\left(-\sec\beta\right)$ & $0.66560453043764058337583145493270$ & $\times$ & $10^{-47}$ \\ [1ex]
 approximation \eqref{eq30} to $a_{25}\left(-\sec\beta\right)$ & $0.66560453043764058337583148788914$ & $\times$ & $10^{-47}$ \\ [1ex]
 error & $-0.3295644$ & $\times$ & $10^{-72}$\\ [1ex] 
 error bound using \eqref{eq33} & $0.6448607$ & $\times$ & $10^{-72}$\\ [-1ex]
 & \\\hline
\end{tabular}
\end{center}
\caption{Approximations for $a_{25}\left(-\sec\beta\right)$ with various $\beta$, using \eqref{eq30}.}
\label{table1}
\end{table*}

More accurate approximations could be derived for the coefficients $a_n\left(-\sec \beta\right)$ by estimating the remainder $A_M \left( n,\beta \right)$ rather than bounding it, but we do not discuss the details here.

\section{Exponentially improved asymptotic expansions}\label{section5}

We shall find it convenient to express our exponentially improved expansions in terms of the (scaled) Terminant function, which is defined by
\[
\widehat T_p \left( w \right) = \frac{{e^{\pi ip} w^{1 - p} e^{ - w} }}{{2\pi i}}\int_0^{ + \infty } {\frac{{t^{p - 1} e^{ - t} }}{w + t}dt} \; \text{ for } \; p>0 \; \text{ and } \; \left| \arg w \right| < \pi ,
\]
and by analytic continuation elsewhere. Olver \cite{Olver4} showed that when $p \sim \left|w\right|$ and $w \to \infty$, we have
\begin{equation}\label{eq36}
ie^{ - \pi ip} \widehat T_p \left( w \right) = \begin{cases} \mathcal{O}\left( {e^{ - w - \left| w \right|} } \right) & \; \text{ if } \; \left| {\arg w} \right| \le \pi \\ \mathcal{O}\left(1\right) & \; \text{ if } \; - 3\pi  < \arg w \le  - \pi. \end{cases}
\end{equation}
Concerning the smooth transition of the Stokes discontinuities, we will use the more precise asymptotics
\begin{equation}\label{eq37}
\widehat T_p \left( w \right) = \frac{1}{2} + \frac{1}{2}\mathop{\text{erf}} \left( {c\left( \varphi  \right)\sqrt {\frac{1}{2}\left| w \right|} } \right) + \mathcal{O}\left( {\frac{{e^{ - \frac{1}{2}\left| w \right|c^2 \left( \varphi  \right)} }}{{\left| w \right|^{\frac{1}{2}} }}} \right)
\end{equation}
for $-\pi +\delta \leq \arg w \leq 3 \pi -\delta$, $0 < \delta  \le 2\pi$; and
\begin{equation}\label{eq38}
e^{ - 2\pi ip} \widehat T_p \left( w \right) =  - \frac{1}{2} + \frac{1}{2}\mathop{\text{erf}} \left( { - \overline {c\left( { - \varphi } \right)} \sqrt {\frac{1}{2}\left| w \right|} } \right) + \mathcal{O}\left( {\frac{{e^{ - \frac{1}{2}\left| w \right|\overline {c^2 \left( { - \varphi } \right)} } }}{{\left| w \right|^{\frac{1}{2}} }}} \right)
\end{equation}
for $- 3\pi  + \delta  \le \arg w \le \pi  - \delta$, $0 < \delta \le 2\pi$. Here $\varphi = \arg w$ and erf denotes the Error function. The quantity $c\left( \varphi  \right)$ is defined implicitly by the equation
\[
\frac{1}{2}c^2 \left( \varphi  \right) = 1 + i\left( {\varphi  - \pi } \right) - e^{i\left( {\varphi  - \pi } \right)},
\]
and corresponds to the branch of $c\left( \varphi  \right)$ which has the following expansion in the neighbourhood of $\varphi = \pi$:
\begin{equation}\label{eq39}
c\left( \varphi  \right) = \left( {\varphi  - \pi } \right) + \frac{i}{6}\left( {\varphi  - \pi } \right)^2  - \frac{1}{{36}}\left( {\varphi  - \pi } \right)^3  - \frac{i}{{270}}\left( {\varphi  - \pi } \right)^4  +  \cdots .
\end{equation}
For complete asymptotic expansions, see Olver \cite{Olver5}. We remark that Olver uses the different notation $F_p \left( w \right) = ie^{ - \pi ip} \widehat T_p \left( w \right)$ for the Terminant function and the other branch of the function $c\left( \varphi  \right)$. For further properties of the Terminant function, see, for example, Paris and Kaminski \cite[Chapter 6]{Paris3}.

\subsection{Proof of the exponentially improved expansions for $\mathbf{A}_{-\nu}\left(\nu x\right)$}

\subsubsection{Case (i): $x>1$} First, we suppose that $\left|\arg \nu\right| < \frac{\pi}{2}$. Our starting point is the representation \eqref{eq8}, written in the form
\begin{equation}\label{eq43}
R_N \left( {\nu ,\beta } \right) = \frac{{\left( { - 1} \right)^N }}{{2\pi \nu ^{2N + 1} }}\int_0^{ + \infty } {\frac{{t^{2N} }}{{1 - it/\nu }}iH_{it}^{\left( 1 \right)} \left( {it\sec \beta } \right)dt} + \frac{{\left( { - 1} \right)^N }}{{2\pi \nu ^{2N + 1} }}\int_0^{ + \infty } {\frac{{t^{2N} }}{{1 + it/\nu }}iH_{it}^{\left( 1 \right)} \left( {it\sec \beta } \right)dt} .
\end{equation}
Let $0 \leq M <2N$ be a fixed integer. We use \eqref{eq40} to expand the function $H_{it}^{\left( 1 \right)} \left( {it\sec \beta } \right)$ under the integrals in \eqref{eq43}, to obtain
\begin{gather}\label{eq44}
\begin{split}
R_N \left( {\nu ,\beta } \right) = \; & i\frac{{e^{ - \frac{\pi }{4}i} }}{{\left( {\frac{1}{2}\nu \pi \tan \beta } \right)^{\frac{1}{2}} }}\sum\limits_{m = 0}^{M - 1} {\left( { - 1} \right)^m \frac{{U_m \left( {i\cot \beta } \right)}}{{\nu ^m }}\left( { - 1} \right)^m \frac{{\left( {i\nu } \right)^{m - 2N - \frac{1}{2}} }}{{2\pi }}\int_0^{ + \infty } {\frac{{t^{2N - m - \frac{1}{2}} e^{ - t\left( {\tan \beta  - \beta } \right)} }}{{1 - it/\nu }}} dt} \\ & - i\frac{{e^{\frac{\pi }{4}i} }}{{\left( {\frac{1}{2}\nu \pi \tan \beta } \right)^{\frac{1}{2}} }}\sum\limits_{m = 0}^{M - 1} {\frac{{U_m \left( {i\cot \beta } \right)}}{{\nu ^m }}\left( { - 1} \right)^m \frac{{\left( { - i\nu } \right)^{m - 2N - \frac{1}{2}} }}{{2\pi }}\int_0^{ + \infty } {\frac{{t^{2N - m - \frac{1}{2}} e^{ - t\left( {\tan \beta  - \beta } \right)} }}{{1 + it/\nu }}} dt} \\ & + R_{N,M} \left( {\nu ,\beta } \right),
\end{split}
\end{gather}
with
\begin{gather}\label{eq45}
\begin{split}
R_{N,M} \left( {\nu ,\beta } \right) = & - \frac{1}{{\left( {\frac{1}{2}\pi \tan \beta } \right)^{\frac{1}{2}} \left( {i\nu } \right)^{2N + 1} }}\frac{1}{{2\pi i}}\int_0^{ + \infty } {\frac{{t^{2N - \frac{1}{2}} e^{ - t\left( {\tan \beta  - \beta } \right)} }}{{1 - it/\nu }}R_M^{\left( H \right)} \left( {it,\beta } \right)dt}\\ & - \frac{1}{{\left( {\frac{1}{2}\pi \tan \beta } \right)^{\frac{1}{2}} \left( {i\nu } \right)^{2N + 1} }}\frac{1}{{2\pi i}}\int_0^{ + \infty } {\frac{{t^{2N - \frac{1}{2}} e^{ - t\left( {\tan \beta  - \beta } \right)} }}{{1 + it/\nu }}R_M^{\left( H \right)} \left( {it,\beta } \right)dt} .
\end{split}
\end{gather}
The integrals in \eqref{eq44} can be identified in terms of the Terminant function since
\[
\left( { - 1} \right)^m \frac{{\left( {i\nu } \right)^{m - 2N - \frac{1}{2}} }}{{2\pi }}\int_0^{ + \infty } {\frac{{t^{2N - m - \frac{1}{2}} e^{ - t\left( {\tan \beta  - \beta } \right)} }}{{1 - it/\nu }}} dt = e^{i\nu \left( {\tan \beta  - \beta } \right)} \widehat T_{2N - m + \frac{1}{2}} \left( {i\nu \left( {\tan \beta  - \beta } \right)} \right)
\]
and
\[
\left( { - 1} \right)^m \frac{{\left( { - i\nu } \right)^{m - 2N - \frac{1}{2}} }}{{2\pi }}\int_0^{ + \infty } {\frac{{t^{2N - m - \frac{1}{2}} e^{ - t\left( {\tan \beta  - \beta } \right)} }}{{1 + it/\nu }}} dt = e^{ - i\nu \left( {\tan \beta  - \beta } \right)} \widehat T_{2N - m + \frac{1}{2}} \left( { - i\nu \left( {\tan \beta  - \beta } \right)} \right) .
\]
Therefore, we have the following expansion
\begin{align*}
R_N \left( {\nu ,\beta } \right) = \; & i\frac{{e^{i\nu \left( {\tan \beta  - \beta } \right) - \frac{\pi }{4}i} }}{{\left( {\frac{1}{2}\nu \pi \tan \beta } \right)^{\frac{1}{2}} }}\sum\limits_{m = 0}^{M - 1} {\left( { - 1} \right)^m \frac{{U_m \left( {i\cot \beta } \right)}}{{\nu ^m }}\widehat T_{2N - m + \frac{1}{2}} \left( {i\nu \left( {\tan \beta  - \beta } \right)} \right)}\\ & - i\frac{{e^{ - i\nu \left( {\tan \beta  - \beta } \right) + \frac{\pi }{4}i} }}{{\left( {\frac{1}{2}\nu \pi \tan \beta } \right)^{\frac{1}{2}} }}\sum\limits_{m = 0}^{M - 1} {\frac{{U_m \left( {i\cot \beta } \right)}}{{\nu ^m }}\widehat T_{2N - m + \frac{1}{2}} \left( { - i\nu \left( {\tan \beta  - \beta } \right)} \right)}  + R_{N,M} \left( {\nu ,\beta } \right).
\end{align*}
Taking $\nu  = re^{i\theta }$, the representation \eqref{eq45} takes the form
\begin{gather}\label{eq46}
\begin{split}
R_{N,M} \left( {\nu ,\beta } \right) = & - \frac{1}{{\left( {\frac{1}{2}r\pi \tan \beta } \right)^{\frac{1}{2}} \left( {ie^{i\theta } } \right)^{2N + 1} }}\frac{1}{{2\pi i}}\int_0^{ + \infty } {\frac{{\tau ^{2N - \frac{1}{2}} e^{ - r\tau \left( {\tan \beta  - \beta } \right)} }}{{1 - i\tau e^{ - i\theta } }}R_M^{\left( H \right)} \left( {ir\tau ,\beta } \right)d\tau } \\
& - \frac{1}{{\left( {\frac{1}{2}r\pi \tan \beta } \right)^{\frac{1}{2}} \left( {ie^{i\theta } } \right)^{2N + 1} }}\frac{1}{{2\pi i}}\int_0^{ + \infty } {\frac{{\tau ^{2N - \frac{1}{2}} e^{ - r\tau \left( {\tan \beta  - \beta } \right)} }}{{1 + i\tau e^{ - i\theta } }}R_M^{\left( H \right)} \left( {ir\tau ,\beta } \right)d\tau } .
\end{split}
\end{gather}
Using the integral formula \eqref{eq42}, $R_M^{\left( H \right)} \left( {ir\tau ,\beta } \right)$ can be written as
\begin{multline*}
R_M^{\left( H \right)} \left( {ir\tau ,\beta } \right) = \frac{{\left( { - 1} \right)^M }}{{2\left( {2\pi \cot \beta } \right)^{\frac{1}{2}}\left( {r\tau } \right)^M }}\left( {\int_0^{ + \infty } {\frac{{s^{M - \frac{1}{2}} e^{ - s\left( {\tan \beta  - \beta } \right)} }}{{1 + s/r}}\left( {1 + e^{ - 2\pi s} } \right)iH_{is}^{\left( 1 \right)} \left( {is\sec \beta } \right)ds}}\right. \\  + \left.{\left( {\tau  - 1} \right)\int_0^{ + \infty } {\frac{{s^{M - \frac{1}{2}} e^{ - s\left( {\tan \beta  - \beta } \right)} }}{{\left( 1+r\tau /s \right)\left( {1 + s/r} \right)}}\left( {1 + e^{ - 2\pi s} } \right)iH_{is}^{\left( 1 \right)} \left( {is\sec \beta } \right)ds} } \right).
\end{multline*}
Noting that
\[
0 < \frac{1}{{1 + s/r}},\frac{1}{{\left( 1 + r\tau /s \right)\left( {1 + s/r} \right)}} < 1
\]
for positive $r$, $\tau$ and $s$, substitution into \eqref{eq46} yields the upper bound
\begin{align*}
\left| {R_{N,M} \left( {\nu ,\beta } \right)} \right| \le \; & \frac{1}{{\left( {\frac{1}{2}\left| \nu  \right|\pi \tan \beta } \right)^{\frac{1}{2}} }}\frac{{\left| {U_M \left( {i\cot \beta } \right)} \right|}}{{\left| \nu  \right|^M }}\left| {\frac{1}{{2\pi }}\int_0^{ + \infty } {\frac{{\tau ^{2N - M - \frac{1}{2}} e^{ - r\tau \left( {\tan \beta  - \beta } \right)} }}{{1 - i\tau e^{ - i\theta } }}d\tau } } \right|\\
& + \frac{1}{{\left( {\frac{1}{2}\left| \nu  \right|\pi \tan \beta } \right)^{\frac{1}{2}} }}\frac{{\left| {U_M \left( {i\cot \beta } \right)} \right|}}{{\left| \nu  \right|^M }}\frac{1}{{2\pi }}\int_0^{ + \infty } {\tau ^{2N - M - \frac{1}{2}} e^{ - r\tau \left( {\tan \beta  - \beta } \right)} \left| {\frac{{\tau  - 1}}{{\tau  + ie^{i\theta } }}} \right|d\tau } \\
& + \frac{1}{{\left( {\frac{1}{2}\left| \nu  \right|\pi \tan \beta } \right)^{\frac{1}{2}} }}\frac{{\left| {U_M \left( {i\cot \beta } \right)} \right|}}{{\left| \nu  \right|^M }}\left| {\frac{1}{{2\pi }}\int_0^{ + \infty } {\frac{{\tau ^{2N - M - \frac{1}{2}} e^{ - r\tau \left( {\tan \beta  - \beta } \right)} }}{{1 + i\tau e^{ - i\theta } }}d\tau } } \right|\\
& + \frac{1}{{\left( {\frac{1}{2}\left| \nu  \right|\pi \tan \beta } \right)^{\frac{1}{2}} }}\frac{{\left| {U_M \left( {i\cot \beta } \right)} \right|}}{{\left| \nu  \right|^M }}\frac{1}{{2\pi }}\int_0^{ + \infty } {\tau ^{2N - M - \frac{1}{2}} e^{ - r\tau \left( {\tan \beta  - \beta } \right)} \left| {\frac{{\tau  - 1}}{{\tau  - ie^{i\theta } }}} \right|d\tau } .
\end{align*}
Since $\left| {\left( {\tau  - 1} \right)/\left( {\tau  \pm ie^{i\theta } } \right)} \right| \le 1$, we find that
\begin{align*}
\left| {R_{N,M} \left( {\nu ,\beta } \right)} \right| \le \; & \frac{1}{{\left( {\frac{1}{2}\left| \nu  \right|\pi \tan \beta } \right)^{\frac{1}{2}} }}\frac{{\left| {U_M \left( {i\cot \beta } \right)} \right|}}{{\left| \nu  \right|^M }}\left| {e^{i\nu \left( {\tan \beta  - \beta } \right)} \widehat T_{2N - M + \frac{1}{2}} \left( {i\nu \left( {\tan \beta  - \beta } \right)} \right)} \right|\\
& + \frac{1}{{\left( {\frac{1}{2}\left| \nu  \right|\pi \tan \beta } \right)^{\frac{1}{2}} }}\frac{{\left| {U_M \left( {i\cot \beta } \right)} \right|}}{{\left| \nu  \right|^M }}\left| {e^{ - i\nu \left( {\tan \beta  - \beta } \right)} \widehat T_{2N - M + \frac{1}{2}} \left( { - i\nu \left( {\tan \beta  - \beta } \right)} \right)} \right|\\
& + \frac{1}{{\left( {\frac{1}{2}\pi \tan \beta } \right)^{\frac{1}{2}} }}\frac{{\left| {U_M \left( {i\cot \beta } \right)} \right|\Gamma \left( {2N - M + \frac{1}{2}} \right)}}{{\pi \left( {\tan \beta  - \beta } \right)^{2N - M + \frac{1}{2}} \left| \nu  \right|^{2N + 1} }}.
\end{align*}
By continuity, this bound holds in the closed sector $\left|\arg \nu\right| \leq \frac{\pi}{2}$. Assume that $
N = \frac{1}{2}\left| \nu  \right|\left( {\tan \beta  - \beta } \right) + \rho$ where $\rho$ is bounded. Employing Stirling's formula, we find that
\[
\frac{1}{{\left( {\frac{1}{2}\pi \tan \beta } \right)^{\frac{1}{2}} }}\frac{{\left| {U_M \left( {i\cot \beta } \right)} \right|\Gamma \left( {2N - M + \frac{1}{2}} \right)}}{{\pi \left( {\tan \beta  - \beta } \right)^{2N - M + \frac{1}{2}} \left| \nu  \right|^{2N + 1} }} = \mathcal{O}_{M,\rho} \left( {\frac{1}{{\left( {\left| \nu  \right|\left( {\tan \beta  - \beta } \right)} \right)^{\frac{1}{2}} }}\frac{{e^{ - \left| \nu  \right|\left( {\tan \beta  - \beta } \right)} }}{{\left( {\frac{1}{2}\left| \nu  \right|\pi \tan \beta } \right)^{\frac{1}{2}} }}\frac{{\left| {U_M \left( {i\cot \beta } \right)} \right|}}{{\left| \nu  \right|^M }}} \right)
\]
as $\nu \to \infty$. Olver's estimation \eqref{eq36} shows that
\[
\left| {e^{ \pm i\nu \left( {\tan \beta  - \beta } \right)} \widehat T_{2N - M + \frac{1}{2}} \left( { \pm i\nu \left( {\tan \beta  - \beta } \right)} \right)} \right| = \mathcal{O}_{M,\rho } \left( {e^{ - \left| \nu  \right|\left( {\tan \beta  - \beta } \right)} } \right)
\]
for large $\nu$. Therefore, we have that
\begin{equation}\label{eq47}
R_{N,M} \left( {\nu ,\beta } \right) = \mathcal{O}_{M,\rho } \left( {\frac{{e^{ - \left| \nu  \right|\left( {\tan \beta  - \beta } \right)} }}{{\left( {\frac{1}{2}\left| \nu  \right|\pi \tan \beta } \right)^{\frac{1}{2}} }}\frac{{\left| {U_M \left( {i\cot \beta } \right)} \right|}}{{\left| \nu  \right|^M }}} \right)
\end{equation}
as $\nu \to \infty$ in the sector $\left|\arg \nu\right| \leq \frac{\pi}{2}$.

Rotating the path of integration in \eqref{eq45} and applying the residue theorem yields
\begin{gather}\label{eq48}
\begin{split}
R_{N,M} \left( {\nu ,\beta } \right)  = \; & i\frac{{e^{i\nu \left( {\tan \beta  - \beta } \right) - \frac{\pi }{4}i} }}{{\left( {\frac{1}{2}\nu \pi \tan \beta } \right)^{\frac{1}{2}} }}R_M^{\left( H \right)} \left( {\nu ,\beta } \right)  - \frac{1}{{\left( {\frac{1}{2}\pi \tan \beta } \right)^{\frac{1}{2}} \left( {i\nu } \right)^{2N + 1} }}\frac{1}{{2\pi i}}\int_0^{ + \infty } {\frac{{t^{2N - \frac{1}{2}} e^{ - t\left( {\tan \beta  - \beta } \right)} }}{{1 - it/\nu }}R_M^{\left( H \right)} \left( {it,\beta } \right)dt} \\
& - \frac{1}{{\left( {\frac{1}{2}\pi \tan \beta } \right)^{\frac{1}{2}} \left( {i\nu } \right)^{2N + 1} }}\frac{1}{{2\pi i}}\int_0^{ + \infty } {\frac{{t^{2N - \frac{1}{2}} e^{ - t\left( {\tan \beta  - \beta } \right)} }}{{1 + it/\nu }}R_M^{\left( H \right)} \left( {it,\beta } \right)dt} \\
 = \; & i\frac{{e^{i\nu \left( {\tan \beta  - \beta } \right) - \frac{\pi }{4}i} }}{{\left( {\frac{1}{2}\nu \pi \tan \beta } \right)^{\frac{1}{2}} }}R_M^{\left( H \right)} \left( {\nu ,\beta } \right) - R_{N,M} \left( {\nu e^{ - \pi i} ,\beta } \right)
\end{split}
\end{gather}
when $\frac{\pi}{2} < \arg \nu <\frac{3\pi}{2}$. It follows that
\[
\left| {R_{N,M} \left( {\nu ,\beta } \right)} \right| \le \frac{{e^{ - \Im \left( \nu  \right)\left( {\tan \beta  - \beta } \right)} }}{{\left( {\frac{1}{2}\left| \nu  \right|\pi \tan \beta } \right)^{\frac{1}{2}} }}\left| {R_M^{\left( H \right)} \left( {\nu ,\beta } \right)} \right| + \left| {R_{N,M} \left( {\nu e^{ - \pi i} ,\beta } \right)} \right|
\]
in the closed sector $\frac{\pi}{2} \leq \arg \nu \leq \frac{3\pi}{2}$, using continuity. It was proved in \cite{Nemes} that $R_M^{\left( H \right)} \left( {\nu ,\beta } \right) = \mathcal{O}_M \left( {\left|U_M \left( {i\cot \beta } \right)\right|\left| \nu  \right|^{ - M} } \right)$ as $\nu \to \infty$ in the closed sector $-\frac{\pi}{2} \leq \arg \nu \leq \frac{3\pi}{2}$, whence, by \eqref{eq47}, we deduce that
\begin{align*}
R_{N,M} \left( {\nu ,\beta } \right) & = \mathcal{O}_M \left( {\frac{{e^{ - \Im \left( \nu  \right)\left( {\tan \beta  - \beta } \right)} }}{{\left( {\frac{1}{2}\left| \nu  \right|\pi \tan \beta } \right)^{\frac{1}{2}} }}\frac{{\left| {U_M \left( {i\cot \beta } \right)} \right|}}{{\left| \nu  \right|^M }}} \right) + \mathcal{O}_{M,\rho } \left( {\frac{{e^{ - \left| \nu  \right|\left( {\tan \beta  - \beta } \right)} }}{{\left( {\frac{1}{2}\left| \nu  \right|\pi \tan \beta } \right)^{\frac{1}{2}} }}\frac{{\left| {U_M \left( {i\cot \beta } \right)} \right|}}{{\left| \nu  \right|^M }}} \right) \\ &= \mathcal{O}_{M,\rho} \left( {\frac{{e^{ - \Im \left( \nu  \right)\left( {\tan \beta  - \beta } \right)} }}{{\left( {\frac{1}{2}\left| \nu  \right|\pi \tan \beta } \right)^{\frac{1}{2}} }}\frac{{\left| {U_M \left( {i\cot \beta } \right)} \right|}}{{\left| \nu  \right|^M }}} \right)
\end{align*}
as $\nu \to \infty$ in the sector $\frac{\pi}{2} \leq \arg \nu \leq \frac{3\pi}{2}$.

The reflection principle gives the relation
\begin{gather}\label{eq49}
\begin{split}
R_{N,M} \left( {\nu ,\beta } \right) = \overline {R_{N,M} \left( {\bar \nu ,\beta } \right)} & =  - i\frac{{e^{ - i\nu \left( {\tan \beta  - \beta } \right) + \frac{\pi }{4}i} }}{{\left( {\frac{1}{2}\nu \pi \tan \beta } \right)^{\frac{1}{2}} }}\overline {R_M^{\left( H \right)} \left( {\bar \nu ,\beta } \right)}  - R_{N,M} \left( {\nu e^{\pi i} ,\beta } \right) \\ & =  - i\frac{{e^{ - i\nu \left( {\tan \beta  - \beta } \right) + \frac{\pi }{4}i} }}{{\left( {\frac{1}{2}\nu \pi \tan \beta } \right)^{\frac{1}{2}} }}R_M^{\left( H \right)} \left( { \nu e^{\pi i} ,\beta } \right) - R_{N,M} \left( {\nu e^{\pi i} ,\beta } \right),
\end{split}
\end{gather}
valid when $-\frac{3\pi}{2} < \arg \nu < -\frac{\pi}{2}$. Trivial estimation and a continuity argument show that
\[
\left| {R_{N,M} \left( {\nu ,\beta } \right)} \right| \le \frac{{e^{\Im \left( \nu  \right)\left( {\tan \beta  - \beta } \right)} }}{{\left( {\frac{1}{2}\left| \nu  \right|\pi \tan \beta } \right)^{\frac{1}{2}} }}\left| {R_M^{\left( H \right)} \left( { \nu e^{\pi i} ,\beta } \right)} \right| + \left| {R_{N,M} \left( {\nu e^{\pi i} ,\beta } \right)} \right|
\]
in the closed sector $-\frac{3\pi}{2} \leq \arg \nu \leq -\frac{\pi}{2}$. Since $R_M^{\left( H \right)} \left( {\nu e^{\pi i},\beta } \right) = \mathcal{O}_M \left( {\left|U_M \left( {i\cot \beta } \right)\right|\left| \nu  \right|^{ - M} } \right)$ as $\nu \to \infty$ in the closed sector $-\frac{3\pi}{2} \leq \arg \nu \leq \frac{\pi}{2}$, by \eqref{eq47}, we find that
\begin{align*}
R_{N,M} \left( {\nu ,\beta } \right) & = \mathcal{O}_M \left( {\frac{{e^{ \Im \left( \nu  \right)\left( {\tan \beta  - \beta } \right)} }}{{\left( {\frac{1}{2}\left| \nu  \right|\pi \tan \beta } \right)^{\frac{1}{2}} }}\frac{{\left| {U_M \left( {i\cot \beta } \right)} \right|}}{{\left| \nu  \right|^M }}} \right) + \mathcal{O}_{M,\rho } \left( {\frac{{e^{ - \left| \nu  \right|\left( {\tan \beta  - \beta } \right)} }}{{\left( {\frac{1}{2}\left| \nu  \right|\pi \tan \beta } \right)^{\frac{1}{2}} }}\frac{{\left| {U_M \left( {i\cot \beta } \right)} \right|}}{{\left| \nu  \right|^M }}} \right) \\ &= \mathcal{O}_{M,\rho} \left( {\frac{{e^{ \Im \left( \nu  \right)\left( {\tan \beta  - \beta } \right)} }}{{\left( {\frac{1}{2}\left| \nu  \right|\pi \tan \beta } \right)^{\frac{1}{2}} }}\frac{{\left| {U_M \left( {i\cot \beta } \right)} \right|}}{{\left| \nu  \right|^M }}} \right)
\end{align*}
as $\nu \to \infty$ with $-\frac{3\pi}{2} \leq \arg \nu \leq -\frac{\pi}{2}$.

\subsubsection{Case (ii): $x=1$} First, we suppose that $\left|\arg \nu\right| < \frac{\pi}{2}$. We write \eqref{eq19} with $N=0$ in the form
\begin{align*}
\mathbf{A}_{ - \nu } \left( \nu  \right) = \frac{1}{{3\pi \nu ^{\frac{1}{3}} }}\int_0^{ + \infty } {\frac{{t^{ - \frac{2}{3}} e^{ - 2\pi t} }}{{1 + \left( {t/\nu } \right)^2 }}iH_{it}^{\left( 1 \right)} \left( {it} \right)dt} & - \frac{1}{{3\pi \nu }}\int_0^{ + \infty } {\frac{{e^{ - 2\pi t} }}{{1 + \left( {t/\nu } \right)^2 }}iH_{it}^{\left( 1 \right)} \left( {it} \right)dt} \\ & + \frac{1}{{3\pi \nu ^{\frac{5}{3}} }}\int_0^{ + \infty } {\frac{{t^{\frac{2}{3}} e^{ - 2\pi t} }}{{1 + \left( {t/\nu } \right)^2 }}iH_{it}^{\left( 1 \right)} \left( {it} \right)dt} .
\end{align*}
Let $N$, $M$ and $K$ be arbitrary positive integers. Using the expression \eqref{eq2}, we find that
\begin{align*}
\mathbf{A}_{ - \nu } \left( \nu  \right) = \frac{1}{{3\pi \nu ^{\frac{1}{3}} }}\sum\limits_{n = 0}^{N - 1} {d_{6n} \frac{{\Gamma \left( {2n + \frac{1}{3}} \right)}}{{\nu ^{2n} }}} & + \frac{1}{{3\pi \nu }}\sum\limits_{m = 0}^{M - 1} {d_{6m + 2} \frac{{\Gamma \left( {2m + 1} \right)}}{{\nu ^{2m} }}} \\ & + \frac{1}{{3\pi \nu ^{\frac{5}{3}} }}\sum\limits_{k = 0}^{K - 1} {d_{6k + 4} \frac{{\Gamma \left( {2k + \frac{5}{3}} \right)}}{{\nu ^{2k} }}}  + R_{N,M,K} \left( \nu  \right),
\end{align*}
where
\begin{gather}\label{eq55}
\begin{split}
R_{N,M,K} \left( \nu  \right) = \frac{{\left( { - 1} \right)^N }}{{3\pi \nu ^{2N + \frac{1}{3}} }}\int_0^{ + \infty } {\frac{{t^{2N - \frac{2}{3}} e^{ - 2\pi t} }}{{1 + \left( {t/\nu } \right)^2 }}iH_{it}^{\left( 1 \right)} \left( {it} \right)dt} & - \frac{{\left( { - 1} \right)^M }}{{3\pi \nu ^{2M + 1} }}\int_0^{ + \infty } {\frac{{t^{2M} e^{ - 2\pi t} }}{{1 + \left( {t/\nu } \right)^2 }}iH_{it}^{\left( 1 \right)} \left( {it} \right)dt} \\ & + \frac{{\left( { - 1} \right)^K }}{{3\pi \nu ^{2K + \frac{5}{3}} }}\int_0^{ + \infty } {\frac{{t^{2K + \frac{2}{3}} e^{ - 2\pi t} }}{{1 + \left( {t/\nu } \right)^2 }}iH_{it}^{\left( 1 \right)} \left( {it} \right)dt} .
\end{split}
\end{gather}
We remark that $R_{N,N,N} \left( \nu  \right) = R_{3N} \left( \nu  \right)$. Assume that $J$, $L$ and $Q$ are integers such that $0 \leq L < 3N$, $0 \leq L < 3M+1$, $0 \leq Q < 3K+2$ and $J,L,Q \equiv 0 \mod 3$. We apply \eqref{eq53} to expand the function $H_{it}^{\left( 1 \right)} \left( {it} \right)$ under the integral in \eqref{eq55}, to obtain
\begin{gather}\label{eq56}
\begin{split}
R_{N,M,K} \left( \nu  \right) = \; & \frac{2}{{9\pi }}\sum\limits_{j = 0}^{J - 1} {d_{2j} \sin \left( {\frac{{\left( {2j + 1} \right)\pi }}{3}} \right)\frac{{\Gamma \left( {\frac{{2j + 1}}{3}} \right)}}{{\nu ^{\frac{{2j + 1}}{3}} }}\left( { - 1} \right)^{N + j} \frac{{\nu ^{\frac{{2j}}{3} - 2N} }}{\pi }\int_0^{ + \infty } {\frac{{t^{2N - \frac{{2j}}{3} - 1} e^{ - 2\pi t} }}{{1 + \left( {t/\nu } \right)^2 }}dt} } \\
& - \frac{2}{{9\pi }}\sum\limits_{\ell  = 0}^{L - 1} {d_{2\ell } \sin \left( {\frac{{\left( {2\ell  + 1} \right)\pi }}{3}} \right)\frac{{\Gamma \left( {\frac{{2\ell  + 1}}{3}} \right)}}{{\nu ^{\frac{{2\ell  + 1}}{3}} }}\left( { - 1} \right)^{M + \ell } \frac{{\nu ^{\frac{{2\ell  - 2}}{3} - 2M} }}{\pi }\int_0^{ + \infty } {\frac{{t^{2M - \frac{{2\ell  - 2}}{3} - 1} e^{ - 2\pi t} }}{{1 + \left( {t/\nu } \right)^2 }}} dt} \\
&+ \frac{2}{{9\pi }}\sum\limits_{q = 0}^{Q - 1} {d_{2q} \sin \left( {\frac{{\left( {2q + 1} \right)\pi }}{3}} \right)\frac{{\Gamma \left( {\frac{{2q + 1}}{3}} \right)}}{{\nu ^{\frac{{2q + 1}}{3}} }}\left( { - 1} \right)^{K + q} \frac{{\nu ^{\frac{{2q - 4}}{3} - 2K} }}{\pi }\int_0^{ + \infty } {\frac{{t^{2K - \frac{{2q - 4}}{3} - 1} e^{ - 2\pi t} }}{{1 + \left( {t/\nu } \right)^2 }}dt} } \\
& + R_{N,M,K}^{J,L,Q} \left( \nu  \right),
\end{split}
\end{gather}
with
\begin{gather}\label{eq57}
\begin{split}
R_{N,M,K}^{J,L,Q} \left( \nu  \right) = \frac{{\left( { - 1} \right)^N }}{{3\pi \nu ^{2N + \frac{1}{3}} }}\int_0^{ + \infty } {\frac{{t^{2N - \frac{2}{3}} e^{ - 2\pi t} }}{{1 + \left( {t/\nu } \right)^2 }}iR_J^{\left( H \right)} \left( {it} \right)dt}  & - \frac{{\left( { - 1} \right)^M }}{{3\pi \nu ^{2M + 1} }}\int_0^{ + \infty } {\frac{{t^{2M} e^{ - 2\pi t} }}{{1 + \left( {t/\nu } \right)^2 }}iR_L^{\left( H \right)} \left( {it} \right)dt} \\ & + \frac{{\left( { - 1} \right)^K }}{{3\pi \nu ^{2K + \frac{5}{3}} }}\int_0^{ + \infty } {\frac{{t^{2K + \frac{2}{3}} e^{ - 2\pi t} }}{{1 + \left( {t/\nu } \right)^2 }}iR_Q^{\left( H \right)} \left( {it} \right)dt} .
\end{split}
\end{gather}
The integrals in \eqref{eq56} can be identified in terms of the Terminant function since
\[
\left( { - 1} \right)^{N + j} \frac{{\nu ^{\frac{{2j}}{3} - 2N} }}{\pi }\int_0^{ + \infty } {\frac{{t^{2N - \frac{{2j}}{3} - 1} e^{ - 2\pi t} }}{{1 + \left( {t/\nu } \right)^2 }}dt}  = ie^{ - 2\pi i\nu } \widehat T_{2N - \frac{{2j}}{3}} \left( { - 2\pi i\nu } \right) - ie^{\frac{\pi }{3}i} e^{2\pi i\nu } e^{\frac{{2\left( {2j + 1} \right)\pi i}}{3}} \widehat T_{2N - \frac{{2j}}{3}} \left( {2\pi i\nu } \right),
\]
\begin{align*}
\left( { - 1} \right)^{M + \ell } \frac{{\nu ^{\frac{{2\ell  - 2}}{3} - 2M} }}{\pi }\int_0^{ + \infty } {\frac{{t^{2M - \frac{{2\ell  - 2}}{3} - 1} e^{ - 2\pi t} }}{{1 + \left( {t/\nu } \right)^2 }}} dt = & - ie^{ - 2\pi i\nu } \widehat T_{2M - \frac{{2\ell  - 2}}{3}} \left( { - 2\pi i\nu } \right)\\ & - ie^{2\pi i\nu } e^{\frac{{2\left( {2\ell  + 1} \right)\pi i}}{3}} \widehat T_{2M - \frac{{2\ell  - 2}}{3}} \left( {2\pi i\nu } \right),
\end{align*}
and
\begin{align*}
\left( { - 1} \right)^{K + q} \frac{{\nu ^{\frac{{2q - 4}}{3} - 2K} }}{\pi }\int_0^{ + \infty } {\frac{{t^{2K - \frac{{2q - 4}}{3} - 1} e^{ - 2\pi t} }}{{1 + \left( {t/\nu } \right)^2 }}dt}  =\; & ie^{ - 2\pi i\nu } \widehat T_{2K - \frac{{2q - 4}}{3}} \left( { - 2\pi i\nu } \right) \\ &- ie^{ - \frac{\pi }{3}i} e^{2\pi i\nu } e^{\frac{{2\left( {2q + 1} \right)\pi i}}{3}} \widehat T_{2K - \frac{{2q - 4}}{3}} \left( {2\pi i\nu } \right).
\end{align*}
Substitution into \eqref{eq56} leads to the expansion \eqref{eq68}. Taking $\nu  = re^{i\theta }$, the representation \eqref{eq57} becomes
\begin{gather}\label{eq58}
\begin{split}
R_{N,M,K}^{J,L,Q} \left( \nu  \right) = \; & \Phi _ +  \left( {N,2N + \frac{1}{3},J} \right) + \Phi _ -  \left( {N,2N + \frac{1}{3},J} \right) - \Phi _ +  \left( {M,2M + 1,L} \right) \\ & - \Phi _ -  \left( {M,2M + 1,L} \right) + \Phi _ +  \left( {K,2K + \frac{5}{3},Q} \right) + \Phi _ -  \left( {K,2K + \frac{5}{3},Q} \right),
\end{split}
\end{gather}
with
\[
\Phi _ \pm  \left( {A,B,C} \right) = \frac{{\left( { - 1} \right)^A }}{{6\pi \left( {e^{i\theta } } \right)^B }}\int_0^{ + \infty } {\frac{{\tau ^{B - 1} e^{ - 2\pi r\tau } }}{{1 \pm i\tau e^{ - i\theta } }}iR_C^{\left( H \right)} \left( {ir\tau } \right)d\tau } .
\]
In \cite[Appendix B]{Nemes} it was shown that
\[
\frac{{1 - \left( {s/r\tau } \right)^{\frac{4}{3}} }}{{1 - \left( {s/r\tau } \right)^2 }} = \frac{{1 - \left( {s/r} \right)^{\frac{4}{3}} }}{{1 - \left( {s/r} \right)^2 }} + \left( {\tau  - 1} \right)f\left( {r,\tau ,s} \right)
\]
for positive $r$, $\tau$ and $s$, with some $f\left(r,\tau ,s\right)$ satisfying $\left|f\left(r,\tau ,s\right)\right| \leq 2$. Using the integral formula \eqref{eq54}, $R_J^{\left( H \right)} \left( {ir\tau } \right)$ can be written as
\begin{align*}
R_J^{\left( H \right)} \left( {ir\tau } \right) = \; & \frac{1}{{\sqrt 3 \pi \left( {r\tau } \right)^{\frac{{2J + 1}}{3}} }}\int_0^{ + \infty } {s^{\frac{{2J - 2}}{3}} e^{ - 2\pi s} \frac{{1 - \left( {s/r\tau } \right)^{\frac{4}{3}} }}{{1 - \left( {s/r\tau } \right)^2 }}H_{is}^{\left( 1 \right)} \left( {is} \right)ds} \\
= \; & \frac{1}{{\sqrt 3 \pi \left( {r\tau } \right)^{\frac{{2J + 1}}{3}} }}\int_0^{ + \infty } {s^{\frac{{2J - 2}}{3}} e^{ - 2\pi s} \frac{{1 - \left( {s/r} \right)^{\frac{4}{3}} }}{{1 - \left( {s/r} \right)^2 }}H_{is}^{\left( 1 \right)} \left( {is} \right)ds} \\
& + \frac{{\tau  - 1}}{{\sqrt 3 \pi \left( {r\tau } \right)^{\frac{{2J + 1}}{3}} }}\int_0^{ + \infty } {s^{\frac{{2J - 2}}{3}} e^{ - 2\pi s} f\left( {r,\tau ,s} \right)H_{is}^{\left( 1 \right)} \left( {is} \right)ds} ,
\end{align*}
and similarly for $R_L^{\left( H \right)} \left( {ir\tau } \right)$ and $R_Q^{\left( H \right)} \left( {ir\tau } \right)$. Noting that
\[
0< \frac{{1 - \left( {s/r} \right)^{\frac{4}{3}} }}{{1 - \left( {s/r} \right)^2 }} < 1
\]
for any positive $r$ and $s$, substitution into \eqref{eq58} yields the upper bound
\begin{align*}
\left| {R_{N,M,K}^{J,L,Q} \left( \nu  \right)} \right| \le \; & \Xi _ +  \left( {2N,2J,\frac{{2J}}{3}} \right) + \Xi _ -  \left( {2N,2J,\frac{{2J}}{3}} \right) + \Xi _ +  \left( {2M,2L,\frac{{2L - 2}}{3}} \right) \\ & + \Xi _ -  \left( {2M,2L,\frac{{2L - 2}}{3}} \right) + \Xi _ +  \left( {2K,2Q,\frac{{2Q - 4}}{3}} \right) + \Xi _ -  \left( {2K,2Q,\frac{{2Q - 4}}{3}} \right),
\end{align*}
with
\[
\Xi _ \pm  \left( {A,B,C} \right) = \frac{{\left| {d_B } \right|\Gamma \left( {\frac{{B + 1}}{3}} \right)}}{{3\sqrt 3 \pi \left| \nu  \right|^{\frac{{B + 1}}{3}} }}\left( {\left| {\frac{1}{{2\pi }}\int_0^{ + \infty } {\frac{{\tau ^{A - C - 1} e^{ - 2\pi r\tau } }}{{1 \pm i\tau e^{ - i\theta } }}d\tau } } \right| + \frac{1}{\pi }\int_0^{ + \infty } {\tau ^{A - C - 1} e^{ - 2\pi r\tau } \left| {\frac{{\tau  - 1}}{{\tau  \mp ie^{i\theta } }}} \right|d\tau } } \right).
\]
As $\left| \left(\tau  - 1\right)/\left(\tau \pm i e^{i\theta}\right)  \right| \le 1$, we find that
\begin{align*}
\left| {R_{N,M,K}^{J,L,Q} \left( \nu  \right)} \right| \le \; & \frac{{\left| {d_{2J} } \right|\Gamma \left( {\frac{{2J + 1}}{3}} \right)}}{{3\sqrt 3 \pi \left| \nu  \right|^{\frac{{2J + 1}}{3}} }}\left| {e^{ - 2\pi i\nu } \widehat T_{2N - \frac{{2J}}{3}} \left( { - 2\pi i\nu } \right)} \right| + \frac{{\left| {d_{2J} } \right|\Gamma \left( {\frac{{2J + 1}}{3}} \right)}}{{3\sqrt 3 \pi \left| \nu  \right|^{\frac{{2J + 1}}{3}} }}\left| {e^{2\pi i\nu } \widehat T_{2N - \frac{{2J}}{3}} \left( {2\pi i\nu } \right)} \right|\\
& + \frac{{2\left| {d_{2J} } \right|\Gamma \left( {\frac{{2J + 1}}{3}} \right)\Gamma \left( {2N - \frac{{2J}}{3}} \right)}}{{3\sqrt 3 \pi ^2 \left( {2\pi } \right)^{2N - \frac{{2J}}{3}} \left| \nu  \right|^{2N + \frac{1}{3}} }} + \frac{{\left| {d_{2L} } \right|\Gamma \left( {\frac{{2L + 1}}{3}} \right)}}{{3\sqrt 3 \pi \left| \nu  \right|^{\frac{{2L + 1}}{3}} }}\left| {e^{ - 2\pi i\nu } \widehat T_{2N - \frac{{2L - 2}}{3}} \left( { - 2\pi i\nu } \right)} \right|\\
& + \frac{{\left| {d_{2L} } \right|\Gamma \left( {\frac{{2L + 1}}{3}} \right)}}{{3\sqrt 3 \pi \left| \nu  \right|^{\frac{{2L + 1}}{3}} }}\left| {e^{2\pi i\nu } \widehat T_{2N - \frac{{2L - 2}}{3}} \left( {2\pi i\nu } \right)} \right| + \frac{{2\left| {d_{2L} } \right|\Gamma \left( {\frac{{2L + 1}}{3}} \right)\Gamma \left( {2M - \frac{{2L - 2}}{3}} \right)}}{{3\sqrt 3 \pi ^2 \left( {2\pi } \right)^{2M - \frac{{2L - 2}}{3}} \left| \nu  \right|^{2M + 1} }}\\
& + \frac{{\left| {d_{2Q} } \right|\Gamma \left( {\frac{{2Q + 1}}{3}} \right)}}{{3\sqrt 3 \pi \left| \nu  \right|^{\frac{{2Q + 1}}{3}} }}\left| {e^{ - 2\pi i\nu } \widehat T_{2N - \frac{{2Q - 4}}{3}} \left( { - 2\pi i\nu } \right)} \right| + \frac{{\left| {d_{2Q} } \right|\Gamma \left( {\frac{{2Q + 1}}{3}} \right)}}{{3\sqrt 3 \pi \left| \nu  \right|^{\frac{{2Q + 1}}{3}} }}\left| {e^{2\pi i\nu } \widehat T_{2N - \frac{{2Q - 4}}{3}} \left( {2\pi i\nu } \right)} \right|\\
& + \frac{{2\left| {d_{2Q} } \right|\Gamma \left( {\frac{{2Q + 1}}{3}} \right)\Gamma \left( {2K - \frac{{2Q - 4}}{3}} \right)}}{{3\sqrt 3 \pi ^2 \left( {2\pi } \right)^{2K - \frac{{2Q - 4}}{3}} \left| \nu  \right|^{2K + \frac{5}{3}} }} .
\end{align*}
By continuity, this bound holds in the closed sector $\left|\arg \nu\right| \le \frac{\pi}{2}$. Suppose that $N= \pi\left| \nu  \right| + \rho $, $M= \pi\left| \nu  \right| + \sigma $ and $K= \pi\left| \nu  \right| + \eta$ where $\rho$, $\sigma$ and $\eta$ are bounded. An application of Stirling's formula shows that
\[
\frac{{2\left| {d_{2J} } \right|\Gamma \left( {\frac{{2J + 1}}{3}} \right)\Gamma \left( {2N - \frac{{2J}}{3}} \right)}}{{3\sqrt 3 \pi ^2 \left( {2\pi } \right)^{2N - \frac{{2J}}{3}} \left| \nu  \right|^{2N + \frac{1}{3}} }} = \mathcal{O}_{J,\rho } \left( {\frac{{e^{ - 2\pi \left| \nu  \right|} }}{{\left| \nu  \right|^{\frac{1}{2}} }}\left| {d_{2J} } \right|\frac{{\Gamma \left( {\frac{{2J + 1}}{3}} \right)}}{{\left| \nu  \right|^{\frac{{2J + 1}}{3}} }}} \right),
\]
\[
\frac{{2\left| {d_{2L} } \right|\Gamma \left( {\frac{{2L + 1}}{3}} \right)\Gamma \left( {2M - \frac{{2L - 2}}{3}} \right)}}{{3\sqrt 3 \pi ^2 \left( {2\pi } \right)^{2M - \frac{{2L - 2}}{3}} \left| \nu  \right|^{2M + 1} }} = \mathcal{O}_{L,\sigma } \left( {\frac{{e^{ - 2\pi \left| \nu  \right|} }}{{\left| \nu  \right|^{\frac{1}{2}} }}\left| {d_{2L} } \right|\frac{{\Gamma \left( {\frac{{2L + 1}}{3}} \right)}}{{\left| \nu  \right|^{\frac{{2L + 1}}{3}} }}} \right),
\]
and
\[
\frac{{2\left| {d_{2Q} } \right|\Gamma \left( {\frac{{2Q + 1}}{3}} \right)\Gamma \left( {2K - \frac{{2Q - 4}}{3}} \right)}}{{3\sqrt 3 \pi ^2 \left( {2\pi } \right)^{2K - \frac{{2Q - 4}}{3}} \left| \nu  \right|^{2K + \frac{5}{3}} }} = \mathcal{O}_{Q,\eta } \left( {\frac{{e^{ - 2\pi \left| \nu  \right|} }}{{\left| \nu  \right|^{\frac{1}{2}} }}\left| {d_{2Q} } \right|\frac{{\Gamma \left( {\frac{{2Q + 1}}{3}} \right)}}{{\left| \nu  \right|^{\frac{{2Q + 1}}{3}} }}} \right)
\]
as $\nu \to \infty$. Using Olver's estimation \eqref{eq36}, we find
\[
\left| {e^{ \pm 2\pi i\nu } \widehat T_{2N - \frac{{2J}}{3}} \left( { \pm 2\pi i\nu } \right)} \right| = \mathcal{O}_{J,\rho } \left( {e^{ - 2\pi \left| \nu  \right|} } \right),
\]
\[
\left| {e^{ \pm 2\pi i\nu } \widehat T_{2N - \frac{{2L - 2}}{3}} \left( { \pm 2\pi i\nu } \right)} \right| = \mathcal{O}_{L,\sigma } \left( {e^{ - 2\pi \left| \nu  \right|} } \right),
\]
and
\[
\left| {e^{ \pm 2\pi i\nu } \widehat T_{2N - \frac{{2Q - 4}}{3}} \left( { \pm 2\pi i\nu } \right)} \right| = \mathcal{O}_{Q,\eta } \left( {e^{ - 2\pi \left| \nu  \right|} } \right)
\]
for large $\nu$. Therefore, we have
\begin{gather}\label{eq59}
\begin{split}
R_{N,M,K}^{J,L,Q} \left( \nu  \right) = \mathcal{O}_{J,\rho } \left( {e^{ - 2\pi \left| \nu  \right|} \left| {d_{2J} } \right|\frac{{\Gamma \left( {\frac{{2J + 1}}{3}} \right)}}{{\left| \nu  \right|^{\frac{{2J + 1}}{3}} }}} \right) & + \mathcal{O}_{L,\sigma } \left( {e^{ - 2\pi \left| \nu  \right|} \left| {d_{2L} } \right|\frac{{\Gamma \left( {\frac{{2L + 1}}{3}} \right)}}{{\left| \nu  \right|^{\frac{{2L + 1}}{3}} }}} \right)\\
& + \mathcal{O}_{Q,\eta } \left( {e^{ - 2\pi \left| \nu  \right|} \left| {d_{2Q} } \right|\frac{{\Gamma \left( {\frac{{2Q + 1}}{3}} \right)}}{{\left| \nu  \right|^{\frac{{2Q + 1}}{3}} }}} \right)
\end{split}
\end{gather}
as $\nu \to \infty$ in the sector $\left|\arg\nu\right| \leq \frac{\pi}{2}$.

Next, we consider the sector $\frac{\pi}{2} < \arg \nu  < \frac{3\pi}{2}$. Rotating the path of integration in \eqref{eq57} and applying the residue theorem gives
\begin{gather}\label{eq61}
\begin{split}
R_{N,M,K}^{J,L,Q} \left( \nu  \right) = \; &  ie^{\frac{\pi }{3}i} \frac{{e^{2\pi i\nu } }}{3}R_J^{\left( H \right)} \left( \nu  \right) +  \frac{\left( { - 1} \right)^N}{{3\pi \nu^{2N + \frac{1}{3}} }}\int_0^{ + \infty } {\frac{{t^{2N - \frac{2}{3}} e^{ - 2\pi t} }}{{1 + \left( {t/\nu e^{ - \pi i} } \right)^2 }}iR_J^{\left( H \right)} \left( {it} \right)dt} \\
& - i\frac{{e^{2\pi i\nu } }}{3}R_L^{\left( H \right)} \left( \nu  \right) - \frac{{\left( { - 1} \right)^M }}{{3\pi \nu^{2M + 1} }}\int_0^{ + \infty } {\frac{{t^{2M} e^{ - 2\pi t} }}{{1 + \left( {t/\nu e^{ - \pi i} } \right)^2 }}iR_L^{\left( H \right)} \left( {it} \right)dt} \\
& + ie^{ - \frac{\pi }{3}i} \frac{{e^{2\pi i\nu } }}{3}R_Q^{\left( H \right)} \left( \nu  \right) + \frac{{\left( { - 1} \right)^K }}{{3\pi \nu^{2K + \frac{5}{3}} }}\int_0^{ + \infty } {\frac{{t^{2K + \frac{2}{3}} e^{ - 2\pi t} }}{{1 + \left( {t/\nu e^{ - \pi i} } \right)^2 }}iR_Q^{\left( H \right)} \left( {it} \right)dt} ,
\end{split}
\end{gather}
for $\frac{\pi}{2} < \arg \nu  < \frac{3\pi}{2}$. It is easy to see that the sum of three integrals has the order of magnitude given in the right-hand side of \eqref{eq59}. It follows that when $J=K=Q$, the bound \eqref{eq59} remains valid in the wider sector $-\frac{\pi}{2} \leq \arg \nu  \leq \frac{3\pi}{2}$. Otherwise, we have
\begin{multline*}
\left| {ie^{\frac{\pi }{3}i} \frac{{e^{2\pi i\nu } }}{3}R_J^{\left( H \right)} \left( \nu  \right) - i\frac{{e^{2\pi i\nu } }}{3}R_L^{\left( H \right)} \left( \nu  \right) + ie^{ - \frac{\pi }{3}i} \frac{{e^{2\pi i\nu } }}{3}R_Q^{\left( H \right)} \left( \nu  \right)} \right|\\
 \le \frac{{e^{ - 2\pi \Im \left( \nu  \right)} }}{3}\left| {R_J^{\left( H \right)} \left( \nu  \right)} \right| + \frac{{e^{ - 2\pi \Im \left( \nu  \right)} }}{3}\left| {R_L^{\left( H \right)} \left( \nu  \right)} \right| + \frac{{e^{ - 2\pi \Im \left( \nu  \right)} }}{3}\left| {R_Q^{\left( H \right)} \left( \nu  \right)} \right| .
\end{multline*}
It was proved in \cite{Nemes} that $R_J^{\left( H \right)} \left( \nu  \right) = \mathcal{O}_J\left( {\left| {d_{2J} } \right|\Gamma \left( {\frac{{2J + 1}}{3}} \right)\left| \nu  \right|^{ - \frac{2J + 1}{3}} } \right)$ as $\nu \to \infty$ in the closed sector $-\frac{\pi}{2} \leq \arg \nu \leq \frac{3\pi}{2}$, whence, by \eqref{eq59}, we deduce that
\begin{gather}\label{eq64}
\begin{split}
R_{N,M,K}^{J,L,Q} \left( \nu  \right) = \mathcal{O}_{J,\rho } \left( {e^{ - 2\pi \Im \left( \nu  \right)} \left| {d_{2J} } \right|\frac{{\Gamma \left( {\frac{{2J + 1}}{3}} \right)}}{{\left| \nu  \right|^{\frac{{2J + 1}}{3}} }}} \right) & + \mathcal{O}_{L,\sigma } \left( {e^{ - 2\pi \Im \left( \nu  \right)} \left| {d_{2L} } \right|\frac{{\Gamma \left( {\frac{{2L + 1}}{3}} \right)}}{{\left| \nu  \right|^{\frac{{2L + 1}}{3}} }}} \right)\\
& + \mathcal{O}_{Q,\eta } \left( {e^{ - 2\pi \Im \left( \nu  \right)} \left| {d_{2Q} } \right|\frac{{\Gamma \left( {\frac{{2Q + 1}}{3}} \right)}}{{\left| \nu  \right|^{\frac{{2Q + 1}}{3}} }}} \right)
\end{split}
\end{gather}
as $\nu \to \infty$ in the sector $\frac{\pi}{2} \leq \arg \nu \leq \frac{3\pi}{2}$.

Similarly, if $J=K=Q$, the bound \eqref{eq59} remains valid in the wider sector $-\frac{3\pi}{2} \leq \arg \nu  \leq \frac{\pi}{2}$; and by the foregoing argument, it is true in the larger sector $-\frac{3\pi}{2} \leq \arg \nu  \leq \frac{3\pi}{2}$. Otherwise, we have
\begin{gather}\label{eq60}
\begin{split}
R_{N,M,K}^{J,L,Q} \left( \nu  \right) = \mathcal{O}_{J,\rho } \left( {e^{ 2\pi \Im \left( \nu  \right)} \left| {d_{2J} } \right|\frac{{\Gamma \left( {\frac{{2J + 1}}{3}} \right)}}{{\left| \nu  \right|^{\frac{{2J + 1}}{3}} }}} \right) & + \mathcal{O}_{L,\sigma } \left( {e^{ 2\pi \Im \left( \nu  \right)} \left| {d_{2L} } \right|\frac{{\Gamma \left( {\frac{{2L + 1}}{3}} \right)}}{{\left| \nu  \right|^{\frac{{2L + 1}}{3}} }}} \right)\\
& + \mathcal{O}_{Q,\eta } \left( {e^{ 2\pi \Im \left( \nu  \right)} \left| {d_{2Q} } \right|\frac{{\Gamma \left( {\frac{{2Q + 1}}{3}} \right)}}{{\left| \nu  \right|^{\frac{{2Q + 1}}{3}} }}} \right)
\end{split}
\end{gather}
for large $\nu$ with $-\frac{3\pi}{2} \leq \arg \nu \leq -\frac{\pi}{2}$.

Consider now the sector $\frac{3\pi}{2} < \arg \nu  < \frac{5\pi}{2}$. Rotation of the path of integration in \eqref{eq61} and application of the residue theorem yields
\begin{align*}
R_{N,M,K}^{J,L,Q} \left( \nu  \right) =\; & ie^{\frac{\pi }{3}i} \frac{{e^{2\pi i\nu } }}{3}R_J^{\left( H \right)} \left( \nu  \right) - i\frac{{e^{2\pi i\nu } }}{3}R_L^{\left( H \right)} \left( \nu  \right) + ie^{ - \frac{\pi }{3}i} \frac{{e^{2\pi i\nu } }}{3}R_Q^{\left( H \right)} \left( \nu  \right)\\
& + i\frac{{e^{ - 2\pi i\nu } }}{3}R_J^{\left( H \right)} \left( {\nu e^{ - \pi i} } \right) + i\frac{{e^{ - 2\pi i\nu } }}{3}R_L^{\left( H \right)} \left( {\nu e^{ - \pi i} } \right) + i\frac{{e^{ - 2\pi i\nu } }}{3}R_Q^{\left( H \right)} \left( {\nu e^{ - \pi i} } \right)\\
& + \frac{{\left( { - 1} \right)^N }}{{3\pi \nu ^{2N + \frac{1}{3}} }}\int_0^{ + \infty } {\frac{{t^{2N - \frac{2}{3}} e^{ - 2\pi t} }}{{1 + \left( {t/\nu e^{ - 2\pi i }} \right)^2 }}iR_J^{\left( H \right)} \left( {it} \right)dt}  - \frac{{\left( { - 1} \right)^M }}{{3\pi \nu ^{2M + 1} }}\int_0^{ + \infty } {\frac{{t^{2M} e^{ - 2\pi t} }}{{1 + \left( {t/\nu e^{ - 2\pi i } } \right)^2 }}iR_L^{\left( H \right)} \left( {it} \right)dt} \\
& + \frac{{\left( { - 1} \right)^K }}{{3\pi \nu ^{2K + \frac{5}{3}} }}\int_0^{ + \infty } {\frac{{t^{2K + \frac{2}{3}} e^{ - 2\pi t} }}{{1 + \left( {t/\nu e^{ - 2\pi i } } \right)^2 }}iR_Q^{\left( H \right)} \left( {it} \right)dt} ,
\end{align*}
for $\frac{3\pi}{2} < \arg \nu  < \frac{5\pi}{2}$. It is easy to see that the sum of three integrals has the order of magnitude given in the right-hand side of \eqref{eq59}. It follows that when $J=K=Q$, the bound \eqref{eq60} holds in the sector $\frac{3\pi}{2} \leq \arg \nu \leq \frac{5\pi}{2}$. Otherwise, we need to bound $R_J^{\left( H \right)} \left( \nu  \right)$, $R_L^{\left( H \right)} \left( \nu  \right)$ and $R_Q^{\left( H \right)} \left( \nu  \right)$. From the connection formula
\[
H_\nu ^{\left( 1 \right)} \left( \nu  \right) =  - H_{\nu e^{ - 2\pi i} }^{\left( 1 \right)} \left( {\nu e^{ - 2\pi i} } \right) - H_{\nu e^{ - 2\pi i} }^{\left( 2 \right)} \left( {\nu e^{ - 2\pi i} } \right) - e^{ - 2\pi i\nu } H_{\nu e^{ - 2\pi i} }^{\left( 2 \right)} \left( {\nu e^{ - 2\pi i} } \right),
\]
we obtain the relation
\begin{align*}
R_J^{\left( H \right)} \left( \nu  \right) & = R_J^{\left( H \right)} \left( {\nu e^{ - 2\pi i} } \right) + R_J^{\left( H \right)} \left( {\nu e^{ - \pi i} } \right) - e^{ - 2\pi i\nu } H_{\nu e^{ - 2\pi i} }^{\left( 2 \right)} \left( {\nu e^{ - 2\pi i} } \right) \\
& = R_J^{\left( H \right)} \left( {\nu e^{ - 2\pi i} } \right) + R_J^{\left( H \right)} \left( {\nu e^{ - \pi i} } \right) + e^{ - 2\pi i\nu } R_0^{\left( H \right)} \left( {\nu e^{ - \pi i} } \right).
\end{align*}
Since $R_J^{\left( H \right)} \left( \nu  \right) = \mathcal{O}_J\left( {\left| {d_{2J} } \right|\Gamma \left( {\frac{{2J + 1}}{3}} \right)\left| \nu  \right|^{ - \frac{2J + 1}{3}} } \right)$ as $\nu \to \infty$ in the sector $-\frac{\pi}{2} \leq \arg \nu \leq \frac{3\pi}{2}$, we infer that
\[
ie^{\frac{\pi }{3}i} \frac{{e^{2\pi i\nu } }}{3}R_J^{\left( H \right)} \left( \nu  \right) = \mathcal{O}_J \left( {e^{ - 2\pi \Im \left( \nu  \right)} \left| {d_{2J} } \right|\frac{{\Gamma \left( {\frac{{2J + 1}}{3}} \right)}}{{\left| \nu  \right|^{\frac{{2J + 1}}{3}} }}} \right) + \mathcal{O}_J \left( {\left| \nu  \right|^{ - \frac{1}{3}} } \right)
\]
for large $\nu$ with $\frac{3\pi}{2} \leq \arg \nu \leq \frac{5\pi}{2}$. A similar estimation holds for the terms involving $R_L^{\left( H \right)} \left( \nu  \right)$ and $R_Q^{\left( H \right)} \left( \nu  \right)$. The sum of the three terms containing $R_J^{\left( H \right)} \left( {\nu e^{ - \pi i} } \right)$, $R_J^{\left( H \right)} \left( {\nu e^{ - \pi i} } \right)$ and $R_Q^{\left( H \right)} \left( {\nu e^{ - \pi i} } \right)$ has the order of magnitude given in the right-hand side of \eqref{eq60}. Therefore, the final result is
\begin{gather}\label{eq70}
\begin{split}
R_{N,M,K}^{J,L,Q} \left( \nu  \right) = \; & \mathcal{O}_{J,\rho } \left( {\cosh \left( {2\pi \Im \left( \nu  \right)} \right)\left| {d_{2J} } \right|\frac{{\Gamma \left( {\frac{{2J + 1}}{3}} \right)}}{{\left| \nu  \right|^{\frac{{2J + 1}}{3}} }}} \right) + \mathcal{O}_{L,\sigma } \left( {\cosh \left( {2\pi \Im \left( \nu  \right)} \right)\left| {d_{2L} } \right|\frac{{\Gamma \left( {\frac{{2L + 1}}{3}} \right)}}{{\left| \nu  \right|^{\frac{{2L + 1}}{3}} }}} \right)\\
& +\mathcal{O}_{Q,\eta } \left( {\cosh \left( {2\pi \Im \left( \nu  \right)} \right)\left| {d_{2Q} } \right|\frac{{\Gamma \left( {\frac{{2Q + 1}}{3}} \right)}}{{\left| \nu  \right|^{\frac{{2Q + 1}}{3}} }}} \right) + \mathcal{O}_{J,L,Q} \left( {\left| \nu  \right|^{ - \frac{1}{3}} } \right)
\end{split}
\end{gather}
as $\nu \to \infty$ in the sector $\frac{3\pi}{2} \leq \arg \nu \leq \frac{5\pi}{2}$.

Similarly, we find that when $J=L=Q$, the estimate \eqref{eq64} holds in the sector $-\frac{5\pi}{2} \leq \arg \nu \leq -\frac{3\pi}{2}$. Otherwise, it can be shown that the estimation \eqref{eq70} is valid in this sector too.

\subsection{Stokes phenomenon and Berry's transition}

\subsubsection{Case (i): $x>1$} We study the Stokes phenomenon related to the asymptotic expansion of $\mathbf{A}_{ - \nu } \left( {\nu \sec \beta } \right)$ occurring when $\arg \nu$ passes through the values $\pm \frac{\pi}{2}$. In the range $\left|\arg \nu\right|<\frac{\pi}{2}$, the asymptotic expansion
\begin{equation}\label{eq52}
\mathbf{A}_{ - \nu } \left( {\nu \sec \beta } \right) \sim  - \frac{1}{\pi }\sum\limits_{n = 0}^\infty  {\frac{{\left( {2n} \right)!a_n \left( { - \sec \beta } \right)}}{{\nu ^{2n + 1} }}} 
\end{equation}
holds as $\nu \to \infty$. From \eqref{eq48} we have
\begin{align*}
\mathbf{A}_{ - \nu } \left( {\nu \sec \beta } \right) = R_{0,0} \left( {\nu ,\beta } \right) & = i\frac{{e^{i\nu \left( {\tan \beta  - \beta } \right) - \frac{\pi }{4}i} }}{{\left( {\frac{1}{2}\nu \pi \tan \beta } \right)^{\frac{1}{2}} }}R_0^{\left( H \right)} \left( {\nu ,\beta } \right) - R_{0,0} \left( {\nu e^{ - \pi i} ,\beta } \right)\\ & = iH_\nu ^{\left( 1 \right)} \left( {\nu \sec \beta } \right) - \mathbf{A}_{\nu} \left( {\nu e^{ - \pi i} \sec \beta } \right)
\end{align*}
when $\frac{\pi}{2} < \arg \nu < \frac{3\pi}{2}$. Similarly, from \eqref{eq49} we find
\[
\mathbf{A}_{ - \nu } \left( {\nu \sec \beta } \right) =  - iH_\nu ^{\left( 2 \right)} \left( {\nu \sec \beta } \right) - \mathbf{A}_{\nu} \left( {\nu e^{\pi i} \sec \beta } \right)
\]
for $-\frac{3\pi}{2} < \arg \nu < -\frac{\pi}{2}$. For the right-hand sides, we can apply the asymptotic expansions of the Hankel functions and the Anger--Weber function to deduce that
\begin{equation}\label{eq75}
\mathbf{A}_{ - \nu } \left( {\nu \sec \beta } \right) \sim i\frac{{e^{i\nu \left( {\tan \beta  - \beta } \right) - \frac{\pi }{4}i} }}{{\left( {\frac{1}{2}\nu \pi \tan \beta } \right)^{\frac{1}{2}} }}\sum\limits_{m = 0}^\infty  {\left( { - 1} \right)^m \frac{{U_m \left( {i\cot \beta } \right)}}{{\nu ^m }}}  - \frac{1}{\pi }\sum\limits_{n = 0}^\infty  {\frac{{\left( {2n} \right)!a_n \left( { - \sec \beta } \right)}}{{\nu ^{2n + 1} }}} 
\end{equation}
as $\nu \to \infty$ in the sector $\frac{\pi}{2} < \arg \nu < \frac{3\pi}{2}$, and
\begin{equation}\label{eq76}
\mathbf{A}_{ - \nu } \left( {\nu \sec \beta } \right) \sim  - i\frac{{e^{ - i\nu \left( {\tan \beta  - \beta } \right) + \frac{\pi }{4}i} }}{{\left( {\frac{1}{2}\nu \pi \tan \beta } \right)^{\frac{1}{2}} }}\sum\limits_{m = 0}^\infty  {\frac{{U_m \left( {i\cot \beta } \right)}}{{\nu ^m }}}  - \frac{1}{\pi }\sum\limits_{n = 0}^\infty  {\frac{{\left( {2n} \right)!a_n \left( { - \sec \beta } \right)}}{{\nu ^{2n + 1} }}} 
\end{equation}
as $\nu \to \infty$ in the sector $-\frac{3\pi}{2} < \arg \nu < -\frac{\pi}{2}$. Therefore, as the line $\arg \nu = \frac{\pi}{2}$ is crossed, the additional series
\begin{equation}\label{eq50}
i\frac{{e^{i\nu \left( {\tan \beta  - \beta } \right) - \frac{\pi }{4}i} }}{{\left( {\frac{1}{2}\nu \pi \tan \beta } \right)^{\frac{1}{2}} }}\sum\limits_{m = 0}^\infty  {\left( { - 1} \right)^m \frac{{U_m \left( {i\cot \beta } \right)}}{\nu ^m}}
\end{equation}
appears in the asymptotic expansion of $\mathbf{A}_{ - \nu } \left( {\nu \sec \beta } \right)$ beside the original one \eqref{eq52}. Similarly, as we pass through the line $\arg \nu = -\frac{\pi}{2}$, the series
\begin{equation}\label{eq51}
- i\frac{{e^{ - i\nu \left( {\tan \beta  - \beta } \right) + \frac{\pi }{4}i} }}{{\left( {\frac{1}{2}\nu \pi \tan \beta } \right)^{\frac{1}{2}} }}\sum\limits_{m = 0}^\infty  {\frac{{U_m \left( {i\cot \beta } \right)}}{\nu ^m }}
\end{equation}
appears in the asymptotic expansion of $\mathbf{A}_{ - \nu } \left( {\nu \sec \beta } \right)$ beside the original series \eqref{eq52}. We have encountered a Stokes phenomenon with Stokes lines $\arg \nu = \pm\frac{\pi}{2}$.

In the important paper \cite{Berry2}, Berry provided a new interpretation of the Stokes phenomenon; he found that assuming optimal truncation, the transition between compound asymptotic expansions is of Error function type, thus yielding a smooth, although very rapid, transition as a Stokes line is crossed.

Using the exponentially improved expansion given in Theorem \ref{thm3}, we show that the asymptotic expansion of $\mathbf{A}_{ - \nu } \left( {\nu \sec \beta } \right)$ exhibits the Berry transition between the two asymptotic series across the Stokes lines $\arg \nu = \pm\frac{\pi}{2}$. More precisely, we shall find that the first few terms of the series in \eqref{eq50} and \eqref{eq51} ``emerge" in a rapid and smooth way as $\arg \nu$ passes through $\frac{\pi}{2}$ and $-\frac{\pi}{2}$, respectively.

From Theorem \ref{thm3}, we conclude that if $N \approx \frac{1}{2}\left| \nu  \right|\left( {\tan \beta  - \beta } \right)$, then for large $\nu$, $ \left|\arg \nu\right| < \pi$, we have
\begin{align*}
\mathbf{A}_{ - \nu } \left( {\nu \sec \beta } \right) \approx & - \frac{1}{\pi }\sum\limits_{n = 0}^{N - 1} {\frac{{\left( {2n} \right)!a_n \left( { - \sec \beta } \right)}}{\nu ^{2n + 1}}} \\ & + i\frac{{e^{i\nu \left( {\tan \beta  - \beta } \right) - \frac{\pi }{4}i} }}{{\left( {\frac{1}{2}\nu \pi \tan \beta } \right)^{\frac{1}{2}} }}\sum\limits_{m = 0} {\left( { - 1} \right)^m \frac{{U_m \left( {i\cot \beta } \right)}}{{\nu ^m }}\widehat T_{2N - m + \frac{1}{2}} \left( {i\nu \left( {\tan \beta  - \beta } \right)} \right)} 
\\ & - i\frac{{e^{ - i\nu \left( {\tan \beta  - \beta } \right) + \frac{\pi }{4}i} }}{{\left( {\frac{1}{2}\nu \pi \tan \beta } \right)^{\frac{1}{2}} }}\sum\limits_{m = 0} {\frac{{U_m \left( {i\cot \beta } \right)}}{{\nu ^m }}\widehat T_{2N - m + \frac{1}{2}} \left( { - i\nu \left( {\tan \beta  - \beta } \right)} \right)} ,
\end{align*}
where $\sum\nolimits_{m = 0}$ means that the sum is restricted to the leading terms of the series.

In the upper half-plane the terms involving $\widehat T_{2N - m + \frac{1}{2}} \left( { - i\nu \left( {\tan \beta  - \beta } \right)} \right)$ are exponentially small, the dominant contribution comes from the terms involving $\widehat T_{2N - m + \frac{1}{2}} \left( {i\nu \left( {\tan \beta  - \beta } \right)} \right)$. Under the above assumption on $N$, from \eqref{eq37} and \eqref{eq39}, the Terminant functions have the asymptotic behaviour
\[
\widehat T_{2N - m + \frac{1}{2}} \left( {i\nu \left( {\tan \beta  - \beta } \right)} \right) \sim \frac{1}{2} + \frac{1}{2} \mathop{\text{erf}} \left( {\left( {\theta  - \frac{\pi }{2}} \right)\sqrt {\frac{1}{2}\left| \nu  \right|\left( {\tan \beta  - \beta } \right)} } \right)
\]
provided that $\arg \nu = \theta$ is close to $\frac{\pi}{2}$, $\nu$ is large and $m$ is small in comparison with $N$. Therefore, when $\theta  < \frac{\pi}{2}$, the Terminant functions are exponentially small; for $\theta  = \frac{\pi }{2}$, they are asymptotically $\frac{1}{2}$ up to an exponentially small error; and when $\theta  >  \frac{\pi}{2}$, the Terminant functions are asymptotic to $1$ with an exponentially small error. Thus, the transition across the Stokes line $\arg \nu = \frac{\pi}{2}$ is effected rapidly and smoothly. Similarly, in the lower half-plane, the dominant contribution is controlled by the terms involving $\widehat T_{2N - m + \frac{1}{2}} \left( { - i\nu \left( {\tan \beta  - \beta } \right)} \right)$. From \eqref{eq38} and \eqref{eq39}, we have
\[
\widehat T_{2N - m + \frac{1}{2}} \left( { - i\nu \left( {\tan \beta  - \beta } \right)} \right) \sim \frac{1}{2} - \frac{1}{2} \mathop{\text{erf}} \left( {\left( {\theta  + \frac{\pi }{2}} \right)\sqrt {\frac{1}{2}\left| \nu  \right|\left( {\tan \beta  - \beta } \right)} } \right)
\]
under the assumptions that $\arg \nu = \theta$ is close to $-\frac{\pi}{2}$, $\nu$ is large and $m$ is small in comparison with $N \approx \frac{1}{2}\left| \nu  \right|\left( {\tan \beta  - \beta } \right)$. Thus, when $\theta  >  - \frac{\pi}{2}$, the Terminant functions are exponentially small; for $\theta  =  -\frac{\pi}{2}$, they are asymptotic to $\frac{1}{2}$ with an exponentially small error; and when $\theta < - \frac{\pi}{2}$, the Terminant functions are asymptotically $1$ up to an exponentially small error. Therefore, the transition through the Stokes line $\arg \nu = -\frac{\pi}{2}$ is carried out rapidly and smoothly.

We remark that from the expansions \eqref{eq75} and \eqref{eq76}, it follows that \eqref{eq52} is an asymptotic expansion of $\mathbf{A}_{-\nu} \left( {\nu \sec \beta } \right)$ in the wider sector $\left|\arg \nu\right| \leq \pi -\delta < \pi$, with any fixed $0 < \delta \leq \pi$.

\subsubsection{Case (ii): $x=1$} The analysis of the Stokes phenomenon for the asymptotic expansion of $\mathbf{A}_{ - \nu } \left( {\nu } \right)$ is similar to the case $x > 1$. In the range $\left|\arg \nu\right| < \frac{3\pi}{2}$, the asymptotic expansion
\begin{equation}\label{eq65}
\mathbf{A}_{ - \nu } \left( \nu  \right) \sim \frac{1}{{3\pi }}\sum\limits_{n = 0}^{\infty} {d_{2n} \frac{{\Gamma \left( {\frac{{2n + 1}}{3}} \right)}}{{\nu ^{\frac{{2n + 1}}{3}} }}}
\end{equation}
holds as $\nu \to \infty$. Employing the continuation formulas stated in Section \ref{section1}, we find that
\[
\mathbf{A}_{ - \nu } \left( \nu  \right) = \mathbf{A}_{ - \nu } \left( {\nu e^{ - 2\pi i} } \right) - iH_{\nu }^{\left( 1 \right)} \left( {\nu e^{ - 2\pi i} } \right) - ie^{ - 2\pi i\nu } H_{\nu  }^{\left( 2 \right)} \left( {\nu e^{ - 2\pi i} } \right)
\]
and
\[
\mathbf{A}_{ - \nu } \left( \nu  \right) = \mathbf{A}_{ - \nu } \left( {\nu e^{2\pi i} } \right) + ie^{2\pi i\nu } H_{\nu }^{\left( 1 \right)} \left( {\nu e^{2\pi i} } \right) + iH_{\nu }^{\left( 2 \right)} \left( {\nu e^{2\pi i} } \right) .
\]
For the right-hand sides, we can apply the asymptotic expansions of the Hankel functions and the Anger--Weber function to deduce that
\begin{equation}\label{eq77}
\mathbf{A}_{ - \nu } \left( \nu  \right) \sim \frac{1}{{3\pi }}\sum\limits_{n = 0}^\infty  {d_{2n} \frac{{\Gamma \left( {\frac{{2n + 1}}{3}} \right)}}{{\nu ^{\frac{{2n + 1}}{3}} }}}  + ie^{ - 2\pi i\nu } \frac{2}{{3\pi }}\sum\limits_{j = 0}^\infty  {d_{2j} \sin \left( {\frac{{\left( {2j + 1} \right)\pi }}{3}} \right)\frac{{\Gamma \left( {\frac{{2j + 1}}{3}} \right)}}{{\nu ^{\frac{{2j + 1}}{3}} }}} 
\end{equation}
as $\nu \to \infty$ in the sector $\frac{3\pi}{2} < \arg \nu < \frac{5\pi}{2}$, and
\begin{equation}\label{eq78}
\mathbf{A}_{ - \nu } \left( \nu  \right) \sim \frac{1}{{3\pi }}\sum\limits_{n = 0}^\infty  {d_{2n} \frac{{\Gamma \left( {\frac{{2n + 1}}{3}} \right)}}{{\nu ^{\frac{{2n + 1}}{3}} }}}  - ie^{2\pi i\nu } \frac{2}{{3\pi }}\sum\limits_{j = 0}^\infty  {d_{2j} \sin \left( {\frac{{\left( {2j + 1} \right)\pi }}{3}} \right)\frac{{\Gamma \left( {\frac{{2j + 1}}{3}} \right)}}{{\nu ^{\frac{{2j + 1}}{3}} }}} 
\end{equation}
as $\nu \to \infty$ in the sector $-\frac{5\pi}{2} < \arg \nu < -\frac{3\pi}{2}$. Therefore, as the line $\arg \nu = \frac{3\pi}{2}$ is crossed, the additional series
\begin{equation}\label{eq66}
ie^{ - 2\pi i\nu } \frac{2}{{3\pi }}\sum\limits_{j = 0}^\infty  {d_{2j} \sin \left( {\frac{{\left( {2j + 1} \right)\pi }}{3}} \right)\frac{{\Gamma \left( {\frac{{2j + 1}}{3}} \right)}}{{\nu ^{\frac{{2j + 1}}{3}} }}} 
\end{equation}
appears in the asymptotic expansion of $\mathbf{A}_{ - \nu } \left( {\nu } \right)$ beside the original one \eqref{eq65}. Similarly, as we pass through the line $\arg \nu = -\frac{3\pi}{2}$, the series
\begin{equation}\label{eq67}
- ie^{2\pi i\nu } \frac{2}{{3\pi }}\sum\limits_{j = 0}^\infty  {d_{2j} \sin \left( {\frac{{\left( {2j + 1} \right)\pi }}{3}} \right)\frac{{\Gamma \left( {\frac{{2j + 1}}{3}} \right)}}{{\nu ^{\frac{{2j + 1}}{3}} }}} 
\end{equation}
appears in the asymptotic expansion of $\mathbf{A}_{ - \nu } \left( \nu \right)$ beside the original series \eqref{eq65}. We have encountered a Stokes phenomenon with Stokes lines $\arg \nu = \pm\frac{3\pi}{2}$. With the aid of the exponentially improved expansion given in Theorem \ref{thm4}, we shall find that the asymptotic series of $\mathbf{A}_{ - \nu } \left( \nu  \right)$ shows the Berry transition property: the two series in \eqref{eq66} and \eqref{eq67}  "emerge" in a rapid and smooth way as the Stokes lines $\arg \nu = \frac{3\pi}{2}$ and $\arg \nu = -\frac{3\pi}{2}$ are crossed.

Let us assume that in \eqref{eq68} $N,M,K \approx \pi \left|\nu\right|$ and $J=L=Q$. When $\pi  < \arg \nu  < 2\pi$, the terms in \eqref{eq68} involving the Terminant functions of the argument $2\pi i \nu$ are exponentially small, and the main contribution comes from the terms involving the Terminant functions of the argument $-2\pi i \nu$. Therefore, from Theorem \ref{thm4}, we deduce that for large $\nu$, $\pi  < \arg \nu  < 2\pi$, we have
\begin{align*}
& \mathbf{A}_{ - \nu } \left( \nu  \right) \approx \frac{1}{{3\pi \nu ^{\frac{1}{3}} }}\sum\limits_{n = 0}^{N - 1} {d_{6n} \frac{{\Gamma \left( {2n + \frac{1}{3}} \right)}}{{\nu ^{2n} }}} + \frac{1}{{3\pi \nu }}\sum\limits_{m = 0}^{M - 1} {d_{6m + 2} \frac{{\Gamma \left( {2m + 1} \right)}}{{\nu ^{2m} }}}  + \frac{1}{{3\pi \nu ^{\frac{5}{3}} }}\sum\limits_{k = 0}^{K - 1} {d_{6k + 4} \frac{{\Gamma \left( {2k + \frac{5}{3}} \right)}}{{\nu ^{2k} }}}  \\
& + ie^{ - 2\pi i\nu } \frac{2}{{3\pi }}\sum\limits_{j = 0} {d_{2j} \sin \left( {\frac{{\left( {2j + 1} \right)\pi }}{3}} \right)\frac{{\Gamma \left( {\frac{{2j + 1}}{3}} \right)}}{{\nu ^{\frac{{2j + 1}}{3}} }}\frac{{\widehat T_{2N - \frac{{2j}}{3}} \left( { - 2\pi i\nu } \right) + \widehat T_{2M - \frac{{2j - 2}}{3}} \left( { - 2\pi i\nu } \right) + \widehat T_{2K - \frac{{2j - 4}}{3}} \left( { - 2\pi i\nu } \right)}}{3}} ,
\end{align*}
where, as before, $\sum\nolimits_{j = 0}$ means that the sum is restricted to the leading terms of the series.

Since $N,M,K \approx \pi \left|\nu\right|$, from \eqref{eq37} and \eqref{eq39}, the averages of the Terminant functions have the asymptotic behaviour
\[
\frac{{\widehat T_{2N - \frac{{2j}}{3}} \left( { - 2\pi i\nu } \right) + \widehat T_{2M - \frac{{2j - 2}}{3}} \left( { - 2\pi i\nu } \right) + \widehat T_{2K - \frac{{2j - 4}}{3}} \left( { - 2\pi i\nu } \right)}}{3} \sim \frac{1}{2} + \frac{1}{2}\mathop{\text{erf}}\left( {\left( {\theta  - \frac{{3\pi }}{2}} \right)\sqrt {\pi \left| \nu  \right|} } \right),
\]
under the conditions that $\arg \nu = \theta$ is close to $\frac{3\pi}{2}$, $\nu$ is large and $j$ is small compared to $N$, $M$ and $K$. Thus, when $\theta < \frac{3\pi}{2}$, the averages of the Terminant functions are exponentially small; for $\theta  =  \frac{3\pi}{2}$, they are asymptotic to $\frac{1}{2}$ with an exponentially small error; and when $\theta >  \frac{3\pi}{2}$, the averages of the Terminant functions are asymptotically $1$ up to an exponentially small error. Thus, the transition through the Stokes line $\arg \nu = \frac{3\pi}{2}$ is carried out rapidly and smoothly.

Similarly, if $N,M,K \approx \pi \left|\nu\right|$ and $J=L=Q$, then for large $\nu$, $-2\pi  < \arg \nu  < -\pi$, we have
\begin{align*}
\mathbf{A}_{ - \nu } \left( \nu  \right) \approx \; & \frac{1}{{3\pi \nu ^{\frac{1}{3}} }}\sum\limits_{n = 0}^{N - 1} {d_{6n} \frac{{\Gamma \left( {2n + \frac{1}{3}} \right)}}{{\nu ^{2n} }}} + \frac{1}{{3\pi \nu }}\sum\limits_{m = 0}^{M - 1} {d_{6m + 2} \frac{{\Gamma \left( {2m + 1} \right)}}{{\nu ^{2m} }}}  + \frac{1}{{3\pi \nu ^{\frac{5}{3}} }}\sum\limits_{k = 0}^{K - 1} {d_{6k + 4} \frac{{\Gamma \left( {2k + \frac{5}{3}} \right)}}{{\nu ^{2k} }}}  \\
&  - ie^{2\pi i\nu } \frac{2}{{3\pi }}\sum\limits_{j = 0} {d_{2j} \sin \left( {\frac{{\left( {2j + 1} \right)\pi }}{3}} \right)\frac{{\Gamma \left( {\frac{{2j + 1}}{3}} \right)}}{{\nu ^{\frac{{2j + 1}}{3}} }}e^{\frac{{2\left( {2j + 1} \right)\pi i}}{3}}} \\ & \times \frac{e^{\frac{\pi }{3}i} \widehat T_{2N - \frac{{2j}}{3}} \left( {2\pi i\nu } \right) - \widehat T_{2M - \frac{{2j - 2}}{3}} \left( { 2\pi i\nu } \right) + e^{ - \frac{\pi }{3}i} \widehat T_{2K - \frac{{2j - 4}}{3}} \left( { 2\pi i\nu } \right)}{3} .
\end{align*}
From \eqref{eq38} and \eqref{eq39}, the averages of the scaled Terminant functions have the asymptotic behaviour
\[
e^{\frac{{2\left( {2j + 1} \right)\pi i}}{3}} \frac{{e^{\frac{\pi }{3}i} \widehat T_{2N - \frac{{2j}}{3}} \left( {2\pi i\nu } \right) - \widehat T_{2M - \frac{{2j - 2}}{3}} \left( {2\pi i\nu } \right) + e^{ - \frac{\pi }{3}i} \widehat T_{2K - \frac{{2j - 4}}{3}} \left( {2\pi i\nu } \right)}}{3} \sim \frac{1}{2} - \frac{1}{2}\mathop{\text{erf}} \left( {\left( {\theta  + \frac{{3\pi }}{2}} \right)\sqrt {\pi \left| \nu  \right|} } \right),
\]
provided that $N,M,K \approx \pi \left|\nu\right|$, $\arg \nu = \theta$ is close to $-\frac{3\pi}{2}$, $\nu$ is large and $j$ is small compared to $N$, $M$ and $K$. Therefore, when $\theta > - \frac{3\pi}{2}$, the averages of the scaled Terminant functions are exponentially small; for $\theta  =  -\frac{3\pi}{2}$, they are asymptotic to $\frac{1}{2}$ up to an exponentially small error; and when $\theta <  -\frac{3\pi}{2}$, the averages of the scaled Terminant functions are asymptotically $1$ with an exponentially small error. Thus, the transition through the Stokes line $\arg \nu = -\frac{3\pi}{2}$ is effected rapidly and smoothly.

We note that from the expansions \eqref{eq77} and \eqref{eq78}, it follows that \eqref{eq65} is an asymptotic series of $\mathbf{A}_{-\nu} \left( \nu \right)$ in the wider range $\left|\arg \nu\right| \leq 2\pi -\delta < 2\pi$, with any fixed $0 < \delta \leq 2\pi$.

\section{Discussion}\label{section6}

In this paper, we have discussed in detail the large order and argument asymptotics of the Anger--Weber function $\mathbf{A}_{-\nu}\left(\nu x\right)$ when $x \geq 1$, using Howls' method. When $0<x<1$, the path $\mathscr{P}\left(0\right)$, defined in \eqref{eq71}, is not the positive real axis, whence the method is not applicable. If we put $x = \mathop{\text{sech}} \alpha$ with a suitable $\alpha > 0$, the large $\nu$ asymptotics of $\mathbf{A}_{-\nu}\left(\nu x\right)$ can be written as
\begin{equation}\label{eq72}
\mathbf{A}_{ - \nu } \left( {\nu \mathop{\text{sech}} \alpha } \right) \sim \sqrt {\frac{2}{{\pi \nu }}} e^{  \nu \left( {\alpha-\tanh \alpha} \right)} \sum\limits_{n = 0}^\infty  {\frac{{\left( {\frac{1}{2}} \right)_n b_n \left( {\mathop{\text{sech}} \alpha } \right)}}{{\nu ^n }}}
\end{equation}
as $\nu \to +\infty$, with $\left( z \right)_n  = \Gamma \left( {z + n} \right)/\Gamma \left( z \right)$ \cite[p. 298]{NIST}. The first few coefficients are given by
\[
b_0 \left( {\mathop{\text{sech}} \alpha } \right) = \frac{1}{{\left( {1 - \mathop{\text{sech}} ^2 \alpha } \right)^{\frac{1}{4}} }},\; b_1 \left( {\mathop{\text{sech}} \alpha } \right) = \frac{{2 + 3\mathop{\text{sech}} ^2 \alpha }}{{12\left( {1 - \mathop{\text{sech}} ^2 \alpha } \right)^{\frac{7}{4}} }},\; b_2 \left( {\mathop{\text{sech}} \alpha } \right) = \frac{{5 + 300\mathop{\text{sech}} ^2 \alpha  + 81\mathop{\text{sech}} ^4 \alpha }}{{864\left( {1 - \mathop{\text{sech}} ^2 \alpha } \right)^{\frac{{13}}{4}} }}.
\]
It is also known that $\mathbf{A}_{ - \nu } \left( {\nu \mathop{\text{sech}} \alpha } \right)$ has the same asymptotic expansion as the Bessel function $-Y_\nu\left(\nu \mathop{\text{sech}} \alpha\right)$, namely
\begin{equation}\label{eq73}
-Y_\nu  \left( {\nu \mathop{\text{sech}}\alpha } \right) \sim  \frac{{e^{ \nu \left( \alpha-\tanh \alpha \right)} }}{{\left( {\frac{1}{2}\pi \nu \tanh \alpha } \right)^{\frac{1}{2}} }}\sum\limits_{n = 0}^\infty  {\left( { - 1} \right)^n \frac{{U_n \left( {\coth \alpha } \right)}}{{\nu ^n }}} \; \text{ as } \; \nu \to +\infty.
\end{equation}
Here $U_n \left( {\coth \alpha } \right) = \left[U_n \left( x \right)\right]_{x = \coth \alpha }$, where $U_n\left(x\right)$ is a polynomial in $x$ of degree $3n$. These polynomials can be generated by the following recurrence
\[
U_n \left( x \right) = \frac{1}{2}x^2 \left( {1 - x^2 } \right)U'_{n-1} \left( x \right) + \frac{1}{8}\int_0^x {\left( {1 - 5t^2 } \right)U_{n-1} \left( t \right)dt} 
\]
for $n \ge 1$ with $U_0\left(x\right) = 1$ (see, e.g., \cite[p. 376]{Olver}, \cite[p. 256]{NIST}). The uniqueness property of asymptotic power series implies that
\begin{align*}
b_n \left( \mathop{\text{sech}}\alpha  \right) & = \left( { - 1} \right)^n \frac{{2^{2n} n!}}{{\left( {2n} \right)!\tanh ^{\frac{1}{2}} \alpha }}U_n \left( {\coth \alpha } \right)
\\ & = \left( { - 1} \right)^n \frac{{2^{2n} n!}}{{\left( {2n} \right)!\left( {1 - \mathop{\text{sech}}^2\alpha } \right)^{\frac{1}{4}} }}U_n \left( {\left( {1 - {\mathop{\text{sech}}}^2 \alpha } \right)^{ - \frac{1}{2}} } \right)
\end{align*}
for any $n\geq 0$. Based on Darboux's method, Dingle \cite[p. 168]{Dingle} gave a formal asymptotic expansion for the coefficients $U_n \left( {\coth \alpha } \right)$ when $n$ is large. His result, in our notation, may be written as
\begin{equation}\label{eq74}
U_n \left( {\coth \alpha } \right) \approx \frac{{\left( { - 1} \right)^n \Gamma \left( n \right)}}{{2\pi \left( {2\left( {\alpha  - \tanh \alpha } \right)} \right)^n }}\sum\limits_{m = 0}^\infty  {\left( {2\left( {\alpha  - \tanh \alpha } \right)} \right)^m U_m \left( {\coth \alpha } \right)\frac{{\Gamma \left( {n - m} \right)}}{{\Gamma \left( n \right)}}} .
\end{equation}
Numerical calculations indicate that this approximation is correct if it is truncated after the first few terms. Using his formal theory of terminants, Dingle gave exponentially improved versions of \eqref{eq72} and \eqref{eq73} \cite[p. 468 and p. 512]{Dingle}.

As far as we know, no rigorous proof of the late coefficient formula \eqref{eq74} nor realistic error bounds for the expansion \eqref{eq72} are available in the literature. Perhaps, these issues can be handled using differential equation methods, but we leave it as a future research topic.

\appendix
\section{}\label{appendixa}

In this appendix we give some formulas for the computation of the coefficients $a_n \left( -\sec \beta  \right)$ appearing in the large $\nu$ asymptotics of $\mathbf{A}_{-\nu}\left(\nu \sec \beta\right)$. It is known that $a_n \left( { - \sec \beta } \right) = \left[ {a_n \left( \lambda  \right)} \right]_{\lambda  =  - \sec \beta } 
$ where $a_n\left(\lambda\right)$ is a rational function of $\lambda \neq -1$. We consider these rational functions. The first few are given by
\[
a_0 \left( \lambda  \right) = \frac{1}{{1 + \lambda }},\; a_1 \left( \lambda  \right) =  - \frac{\lambda }{{2\left( {1 + \lambda } \right)^4 }},\; a_2 \left( \lambda  \right) = \frac{{9\lambda ^2  - \lambda }}{{24\left( {1 + \lambda } \right)^7 }},\; a_3 \left( \lambda  \right) =  - \frac{{225\lambda ^3  - 54\lambda ^2  + \lambda }}{{720\left( {1 + \lambda } \right)^{10} }} .
\]
From \eqref{eq9} we infer that
\[
a_n \left( \lambda  \right) = \frac{1}{{\left( {2n} \right)!}}\left[ {\frac{{d^{2n} }}{{dt^{2n} }}\left( {\frac{t}{{\lambda \sinh t + t}}} \right)^{2n + 1} } \right]_{t = 0} .
\]

Meijer \cite{Meijer} proved the following explicit formula
\begin{equation}\label{eq34}
a_n \left( \lambda  \right) = \frac{1}{{\left( {1+ \lambda } \right)^{2n + 1} }}\sum\limits_{k = 0}^n {\binom{2n + k}{k}\frac{{\left( { - 1} \right)^k }}{{\left( {2n - 2k} \right)!}}\left[ {\frac{{d^{2n - 2k} }}{{dt^{2n - 2k} }}\left( {\frac{{\sinh t - t}}{{t^3 }}} \right)^k } \right]_{t = 0} \left( {\frac{\lambda }{{1+ \lambda }}} \right)^k } .
\end{equation}
We show that the higher derivatives can be written in terms of the generalized Bernoulli polynomials $B_n^{\left( \kappa  \right)} \left(\ell\right)$, which are defined by the exponential generating function
\[
\left( \frac{z}{e^z  - 1} \right)^\kappa  e^{\ell z}  = \sum\limits_{n = 0}^\infty  {B_n^{\left( \kappa  \right)} \left(\ell\right)\frac{z^n}{n!}} \; \text{ for } \; \left|z\right| < 2\pi.
\]
For basic properties of these polynomials, see Milne-Thomson \cite{Milne-Thomson} or N\"{o}rlund \cite{Norlund}. A straightforward computation gives
\begin{align*}
& \frac{1}{{\left( {2n-2k} \right)!}}\left[ {\frac{{d^{2n-2k} }}{{dt^{2n-2k} }}\left( {\frac{{\sinh t - t}}{{t^3 }}} \right)^k } \right]_{t = 0} = \frac{1}{{2\pi i}}\oint_{\left( {0^ +  } \right)} {\left( {\frac{{\sinh z - z}}{{z^3 }}} \right)^k \frac{{dz}}{{z^{2n-2k + 1} }}} \\
& = \frac{1}{{2\pi i}}\oint_{\left( {0^ +  } \right)} {\left( {\sum\limits_{j = 0}^k {\left( { - 1} \right)^{k - j} \binom{k}{j}z^{k - j} \sinh ^j z} } \right)\frac{{dz}}{{z^{2n + k + 1} }}}  = \sum\limits_{j = 0}^k {\left( { - 1} \right)^{k - j} \binom{k}{j}\frac{1}{{2\pi i}}\oint_{\left( {0^ +  } \right)} {\left( {\frac{{\sinh z}}{z}} \right)^j \frac{dz}{{z^{2n  + 1} }}} } \\
& = \sum\limits_{j = 0}^k {\left( { - 1} \right)^{k - j} \binom{k}{j}\frac{1}{{2\pi i}}\oint_{\left( {0^ +  } \right)} {\left( {\frac{{2z}}{{e^{2z}  - 1}}} \right)^{ - j} e^{ - jz} \frac{dz}{{z^{2n  + 1} }}} }  = \frac{{2^{2n} }}{{\left( {2n} \right)!}}\sum\limits_{j = 0}^k {\left( { - 1} \right)^{k - j} \binom{k}{j}B_{2n}^{\left( { - j} \right)} \left( { - \frac{j}{2}} \right)}.
\end{align*}
Substitution into \eqref{eq34} yields the explicit representation
\[
a_n \left( \lambda  \right) = \frac{{2^{2n} }}{{\left( {2n} \right)!^2 \left( {1+ \lambda} \right)^{2n + 1} }}\sum\limits_{k = 0}^n {\sum\limits_{j = 0}^k {\left( { - 1} \right)^j \frac{{\left( {2n + k} \right)!}}{{\left( {k - j} \right)!j!}}B_{2n}^{\left( { - j} \right)} \left( { - \frac{j}{2}} \right)} \left( {\frac{\lambda }{1+\lambda}} \right)^k } .
\]

In 1952, Lauwerier \cite{Lauwerier} showed that the coefficients in asymptotic expansions of Laplace-type integrals can be calculated by means of linear recurrence relations. Simple application of his method provides the formula
\[
a_n \left( \lambda  \right) = \frac{1}{{\left( {2n} \right)!}}\int_0^{ + \infty } {t^{2n} e^{ - \left( {1 + \lambda } \right)t} P_n \left( t \right)dt} ,
\]
where the polynomials $P_0 \left( x \right), P_1 \left( x \right), P_2 \left( x \right),\ldots$ are given by the recurrence relation
\[
P_n \left( x \right) =  - \sum\limits_{k = 1}^n {\frac{\lambda }{{\left( {2k + 1} \right)!}}\int_0^x {P_{n - k} \left( t \right)dt} } 
\]
with $P_0\left(x\right)=1$.

A simpler recurrence for the $a_n \left( \lambda  \right)$'s can be found using the inhomogeneous Bessel differential equation
\[
\frac{{d^2 w\left( z \right)}}{{dz^2 }} + \frac{1}{z}\frac{{dw  \left( z \right)}}{{dz}} + \left( {1 - \frac{{\nu ^2 }}{{z^2 }}} \right)w\left( z \right) = \frac{{z - \nu }}{{\pi z^2 }}
\]
satisfied by the Anger--Weber function $\mathbf{A}_{\nu}\left(z\right)$. Substituting $z=\nu \lambda$ shows that
\begin{equation}\label{eq35}
\frac{\lambda ^2}{\nu ^2}\frac{{d^2 {\bf A}_\nu  \left( { \nu \lambda} \right)}}{{d\lambda ^2 }} + \frac{\lambda}{\nu ^2}\frac{{d{\bf A}_\nu  \left( { \nu \lambda} \right)}}{d\lambda} + \left( {\lambda ^2  - 1} \right){\bf A}_\nu  \left( {\nu \lambda } \right) = \frac{{\lambda  - 1}}{{\pi \nu }} .
\end{equation}
It is known that for any $\lambda>0$, the function ${\bf A}_\nu  \left( {\nu \lambda} \right)$ has the asymptotic expansion
\[
{\bf A}_\nu  \left( {\nu \lambda} \right) \sim \frac{1}{\pi }\sum\limits_{n = 0}^\infty  {\frac{{\left( {2n} \right)!a_n \left( \lambda  \right)}}{\nu ^{2n + 1}}} 
\]
as $\nu \to \infty$, $\left|\arg \nu\right| <\pi$ (see, e.g, \cite[p. 298]{NIST}). Substituting this series into \eqref{eq35} and equating the coefficients of the inverse powers of $\nu$ we find
\[
a_0 \left( \lambda  \right) = \frac{1}{{1 + \lambda }}\; \text{ and } \;a_n \left( \lambda  \right) = \frac{\lambda}{1 - \lambda ^2} \frac{{\lambda a''_{n - 1} \left( \lambda  \right) + a'_{n - 1} \left( \lambda  \right)}}{2n\left( {2n - 1} \right)} \; \text{ for } \; n\geq 1.
\]

\section*{Acknowledgement} I would like to thank the two anonymous referees for their constructive and helpful comments and suggestions on the manuscript.

\end{document}